\documentclass[english]{smfart}

\usepackage[english]{babel}
 \usepackage{amsfonts}
\usepackage{epsfig,graphics}
\usepackage{amssymb}
\usepackage{amscd}
\usepackage{bull} 

\author{Charles Frances}
\address{IRMA, 7 rue Ren\'e Descartes, 67000 Strasbourg.}
\email{cfrances@math.unistra.fr}
\urladdr{}
\title[Variations on Gromov's open-dense orbit theorem]{Variations on Gromov's open-dense orbit theorem}
\newtheorem{theoreme}{Theorem}[section]

\newtheorem{fact}[theoreme]{Fact}

\newtheorem{proposition}[theoreme]{Proposition}
       
\newtheorem{definition}[theoreme]{Definition}

\newtheorem{lemme}[theoreme]{Lemma}

\newtheorem{question}[theoreme]{Question}

\newtheorem{remarque}[theoreme]{Remark}

\newtheorem{corollaire}[theoreme]{Corollary}

\newtheorem{thmx}{Theorem}

\newcommand{\Ker}{\operatorname{Ker }}

\newcommand{\Hom}{\operatorname{Hom }}

\newcommand{\rk}{\operatorname{rk}}

\newcommand{\ad}{\operatorname{ad}}

\newcommand{\SL}{\operatorname{SL}}

\newcommand{\GL}{\operatorname{GL}}
\newcommand{\PSL}{\operatorname{PSL}}

\newcommand{\SO}{\operatorname{SO}}
\newcommand{\OO}{\operatorname{O}}

\newcommand{\Iso}{\operatorname{Iso}}

\newcommand{\Ad}{\operatorname{Ad}}

\newcommand{\hn}{\hat{N}}
\newcommand{\mathcalm}{\mathcal{M}}

\def\dk{{\mathcal{D}\kappa}}

\def\NN{\mathbb{N}}

\def\RR{\mathbb{R}}

\def\O{\mathcal{O}}
\def\OO{\operatorname{O}}

\newcommand{\mint}{{M^{\rm int}}}

\newcommand{\R}{{\bf R}}

\newcommand{\hm}{{\hat{M}}}
\newcommand{\hx}{{\hat{x}}}
\newcommand{\hy}{{\hat{y}}}

\newcommand{\hgamma}{{\hat{\gamma }}}

\newcommand{\ka}{{\kappa}}

\newcommand{\kil}{\operatorname{{\mathfrak{kill}}}}
\newcommand{\kiloc}{\operatorname{{\mathfrak{kill}^{loc}}}}
\newcommand{\heis}{{\operatorname{\mathfrak{heis}}}}

\newcommand{\lieg}{{\mathfrak{g}}}
\newcommand{\lieh}{{\mathfrak{h}}}
\newcommand{\liep}{{\mathfrak{p}}}

\newcommand{\oo}{{\mathfrak{o}}}

\newcommand{\sld}{\operatorname{{\mathfrak{sl}}(2,\RR)}}

\newenvironment{preuve}{\medskip \noindent {\bf Proof: }}
   {$\diamondsuit$ }

\begin{document}
\frontmatter
\maketitle
\section{Introduction}

The main motivation of this article comes from the following result of M. Gromov, often quoted in the litterature as the
{\it open-dense orbit theorem}.
\begin{theoreme}[\cite{gromov}, Th. 3.3.A]
 \label{thm.gromov}
 Let $M$ be a smooth manifold, and ${\mathcal S}$  a smooth rigid geometric structure of algebraic type on  $M$.  If the automorphism group of 
 $(M,{\mathcal S})$ has a dense orbit, then the structure ${\mathcal S}$ is locally homogeneous on a dense open subset of $M$.
\end{theoreme}

Recall that a structure is locally homogeneous if given any pair of points $(x,y)$ of $M$, there exists a local isometry $f$ 
(namely a local diffeomorphism preserving $\mathcal S$), defined  from a neighborhood of $x$
 to a neighborhood of $y$, and satisfying  $f(x)=y$.  
 
 The  notion of rigid geometric structure of algebraic type was introduced in \cite{gromov}. It covers a wide range of structures, important examples
  of which are 
  pseudo-Riemannian metrics on
  manifolds, or affine connections.  
  
  Theorem \ref{thm.gromov}  actually holds under the weaker assumption that the pseudo-group
  of local isometries
   has a dense orbit.  It comes as a corollary of a more general result, also proved in \cite{gromov}, stating that for a rigid
   geometric structure of algebraic type, there exists a dense open subset where the orbits of the pseudogroup of local isometries are closed submanifolds.  It is thus
    clear that whenever one of these orbits is dense, it must be open. The reader wanting to learn more about Gromov's theory of rigid transformation groups will find
   details in \cite{gromov}, \cite{dambragro}, \cite{benoist}, \cite{zeghibgromov}.  
   
   A beautiful application of Theorem \ref{thm.gromov} can be found in \cite{bfl}, where the authors use Gromov's result to get a full classification
    of contact Anosov flows on compact manifolds, 
   admitting smooth stable and unstable distributions.  Their strategy is to show that such contact flows preserve a smooth 
   pseudo-Riemannian metric, which turns out to be locally homogeneous because of Gromov's theorem, and the Anosov dynamics.  The end 
   (actually the main part) of the proof
    consists in classifying the possible algebraic local models.
More generally, Theorem \ref{thm.gromov} seems to be a key ingredient in   classifiying rigid geometric 
   structures of a certain type, with a topologically
    transtive group of isometries.  There is however a restriction  : the local homogeneity ensured by the theorem is
    only available on a dense 
    open subset of the manifold, while we would like such a result on the whole manifold.  This raises the following natural question :
    
    \begin{question}
     \label{quest.gromov}
     Can the maximal open set of local homogeneity given by Theorem \ref{thm.gromov} be a strict open subset of $M$?
    \end{question}

    While it is expected that the answer to the previous question should be negative, there are very few instances 
     where one can prove it (see   \cite{bfl} and  \cite{dumidolpho} 
     for nontrivial examples where the authors show local homogeneity everywhere). 
     
     \subsection{Open-dense orbit theorem and $3$-dimensional Lorentz metrics}
     
     If we restrict our attention to pseudo-Riemannian structures, the situation seems to be the following. The answer to Question \ref{quest.gromov}  is negative
      for Riemannian manifolds, and pseudo-Riemannian surfaces.  In the first case, it is almost obvious, and in the second one, we see that
       whenever the isometry group has a dense orbit, the sectional curvature must be constant, 
      which implies local homogeneity. 

    One aim of this paper is to study the first nontrivial   case beside the two preceding ones, namely
     that of   $3$-dimensional Lorentz manifolds.  Our main result is :
\begin{thmx}
 \label{thm.orbitedense}
 Let $(M^3,g)$ be a smooth closed $3$-dimensional Lorentz manifold.  
 If the isometry group $\Iso(M,g)$ has a dense orbit, then $(M^3,g)$ is locally homogeneous.
\end{thmx}

Observe that  we don't make any {\it a priori} assumption on the group $\Iso(M,g)$.  In particular, it might be a discrete group. 

Under stronger assumptions, Theorem \ref{thm.orbitedense} can be deduced from previous works.  For instance, if we assume that the metric $g$
 is real analytic, then S. Dumitrescu showed in \cite{sorin.dim3} that the existence
 of a nonempty open orbit for the pseudo-group of local isometries led to local homogeneity.

 
 In the smooth category, under the stronger assumption that there exists a $1$-parameter flow of isometries with
 a dense orbit, Theorem \ref{thm.orbitedense} can be derived from \cite{zeghibflot}, where A. Zeghib classifies completely all Lorentzian flows
  on compact $3$-manifolds which are not equicontinuous.
  
 Let us finally mention that obtaining a generalization of Theorem \ref{thm.orbitedense}  to Lorentz manifolds of arbitrary dimension, or to general 
  pseudo-Riemannian structures seems to be rather challenging. Good examples of topologically transitive pseudo-Riemannian flows, illustrating the 
  conclusions
   of Theorem \ref{thm.gromov}, can be built as follows.
    Let $G$ be a noncompact simple (or semi-simple) Lie group, and let $\Gamma$ be a uniform lattice.  Let $\{ g^t \}$ be a $1$-parameter subgroup 
    of  $G$, whith noncompact closure in $G$.  It follows from Moore's theorem that $\{ g^t \}$ acts ergodically on $G/\Gamma$. Let $\kappa_0$ be the 
    Killing form on $\lieg$.  This is a pseudo-Riemannian scalar product, which is $\Ad (g^t)$-invariant.  Pushing this scalar product by right
     translations, one gets a bi-invariant pseudo-Riemannian metric on $G$, which in turns induces a $g^t$-invariant metric $h_0$ on $G/\Gamma$.
     this Killing metric $h_0$ is actually $G$-invariant, hence homogeneous.  The point is that for many $1$-parameter groups $\{ g^t \}$ (for instance
      when $\{ g^t \}$ is in a Cartan subgroup $A$, or when $\{g^t\}$ is unipotent), there are a lot of pseudo-Riemannian scalar products 
      $\kappa$
       on $\lieg$ beside the Killing form.  Actually, one can choose some $\kappa$'s which are $\Ad(g^t)-$invariant, 
       without being $\Ad(H)-$invariant for $\{g^t\} \subsetneq H$.  By the same construction as above, 
       $\kappa$ yields a pseudo-Riemannian metric $h$ on $G/\Gamma$ for which the isometry group reduces to $\{ g^t \}$
       (maybe up to finite index).  Such a metric is of course no longer homogeneous, but is still locally homogeneous.  
    Now even for those concrete examples, and besides particular choices of groups $\{ g^t \}$, it does not seem obvious to show that {\it all} pseudo-Riemannian
     metrics on $G/\Gamma$ which are $\{ g^t \}-$invariant are locally homogeneous 
     (a generalization of the results of \cite{iozzi} might be usefull to this regard).
%
%


\subsection{Quasihomogeneity}

Let us now discuss purely local problems related to Question \ref{quest.gromov}. 

We recall that a local Killing field for a geometric structure ${\mathcal S}$ on a manifold $M$ is a vector field $X$ defined on some open set
$U \subset M$, and such that the local flow $\varphi_X^t$ preserves ${\mathcal S}$. The set of Killing fields defined on $U$ is a Lie subalgebra
of the vector fields on $U$, that we denote $\kil(U)$.  When the structure $\mathcal S$ is rigid, then $\kil(U)$ is always finite dimensional.
    It then follows that if  $(U_i)_{i \in \NN}$ is a nested family of open sets containing a point $x$, and 
    satisfying $\bigcap_{i \in \NN}U_i=\{x\}$, then the  dimension of $\kil(U_i)$ stabilizes for $i$ large enough.  The resulting Lie algebra
     will be denoted by $\kil(x)$.  
     
Starting from a point $x \in M$, we can consider the set of all points $y \in M$ that can be reached from $x$ by flowing along successive 
 local Killing fields.  This set is called the Kill$^{\rm loc}$-orbit of $x$, and denoted  ${\mathcal O}_x^{\rm loc}$.  
  Let us recall that for ``generic'' rigid structures, there are no local Killing fields at all, and the Kill$^{\rm loc}$-orbits are reduced
   to points. The opposite situation is that of connected  locally homogeneous 
 structures,   for which ${\mathcal O}_x^{\rm loc}=M$.  An interesting weaker notion is   that of quasihomogeneous structure.

\begin{definition} 
A geometric structure is called quasihomogeneous when the union of open Kill$^{\rm loc}$-orbits is dense.  
\end{definition}

 Gromov's theorem \ref{thm.gromov} says that a rigid geometric structure of algebraic type,
with a topologically transitive  automophism group, is quasihomogeneous. It is thus a question of general interest to understand 
when a quasihomogeneous structure
 is actually locally homogeneous.

 It seems that there is no universal  answer to this problem.  For instance, A. Guillot and S. Dumitrescu exhibited
  in \cite{dumidolpho}  quasihomogeneous affine connections on surfaces which are not homogeneous, {\it even in the real analytic category}.
  
  On the contrary,  S. Dumitrescu and K. Melnick recently showed in \cite{sorin.karin} that any real analytic Lorentz metric on a 
  $3$-manifold which is quasihomogeneous must be locally homogeneous.  The analyticity assumption is crucial in their proof, and it is 
    unknown if the results of \cite{sorin.karin} still hold in the smooth category.
    
    Actually, Theorem \ref{thm.orbitedense} will follow from a partial generalization of \cite{sorin.karin} to smooth manifolds. We will 
    indeed show
     the following local result :

\begin{thmx}
 
 \label{thm.quasihomogene}
 Let $(M^3,g)$ be a smooth $3$-dimensional Lorentz manifold (not necessarily closed).  Assume that on a dense open subset, 
 the Lie algebra of local Killing fields
  is at least $4$-dimensional.  Then $(M^3,g)$ is locally homogeneous.  
\end{thmx}

It is not hard to see that the hypothesis on the dimension of the local Killing algebras does imply quasihomogeneity of the metric
(see Fact \ref{fact.quasihomo}). For a quasihomogeneous Lorentz $3$-manifold, the possible dimensions of the local Killing algebras $\kil(x)$ are $3$, 
$4$ or $6$ (in the smooth case, this dimension may vary with the point $x$).  Hence Theorem \ref{thm.quasihomogene} deals with quasihomogeneous 
structures without open Kill$^{\rm loc}$-orbits having a $3$-dimensional local Killing algebra.

 Even if the conclusions of Theorem \ref{thm.quasihomogene} are  the same as for analytic metrics,
the result can not
be obtained by adapting the methods of \cite{sorin.karin}.  Actually, we would like to point out that those regularity issues concentrate a great 
part
 of the subtilities in this kind of problems.  To emphasize this aspect, we observe that the proof of Theorem \ref{thm.quasihomogene} works
  for metrics of class $C^9$ (this regularity is required since we will need several covariant derivatives of the curvature tensor). This 
  regularity is probably not optimal, but let us stress that the  conclusions change 
  dramatically  if we work with metrics which have too low regularity.

 \begin{thmx}
  \label{theo.C}
  There exist $3$-dimensional Lorentz manifolds $(M^3,g)$, such that $g$ is $C^1$ and quasihomogeneous, satisfies hypotheses of 
 Theorem \ref{thm.quasihomogene}, but is not locally homogeneous.
  One can  buildt moreover compact  examples  in regularity $C^0$. 
 \end{thmx}

 The conclusions of Theorem \ref{theo.C}  will be made more precise in section \ref{sec.c1}  (see Theorems \ref{theo.c1} and \ref{thm.quasicompact}).

\subsection{Organization of the paper} 
 A key ingredient in Gromov's theory of rigid transformation groups is a theorem about integration of finite order Killing fields.  We will make
  a systematic use of this result in all our proofs, so that Section \ref{sec.integration} will be devoted to a presentation of this theorem, in the
  convenient framework of Cartan geometries (which includes of course the case of Lorentz metrics).  Next, we will use this integration result
  in Section \ref{sec.criteres}, to prove two general
  criteria allowing to show that
   some quasihomogeneous pseudo-Riemannian structures are actually locally homogeneous. Section 
   \ref{sec.proofquasi} begins with a general study of local Lorentz actions of $4$-dimensional Lie algebras on $3$-manifolds.  This study, together with
    the criteria established in Section \ref{sec.criteres} lead to a proof of  Theorem \ref{thm.quasihomogene}. 
    In Section \ref{sec.dense.orbit}, we  explain how Theorem \ref{thm.orbitedense} can be deduced from
    Theorem \ref{thm.quasihomogene}.  The upshot is to show that when the isometry group of a pseudo-Riemannian manifold is topologically 
    transitive, then numerous local Killing fields must appear (even if the isometry group is discrete, for instance).  Finally, Section \ref{sec.c1} will be
     devoted to the construction of examples of Theorem \ref{theo.C}.


\section{Integration of finite order Killing fields}
\label{sec.integration}

The main tool to understand the  Kill$^{\rm loc}$-orbits of a rigid geometric structure is a theorem about integration of
 finite order Killing fields proved in \cite{gromov}[Section 1.6].  The results of \cite{gromov} generalize former integrability theorems
 proved by K. Nomizu in \cite{nomizu} and I. Singer in \cite{singer}.  We won't follow here the approach of \cite{gromov}, but rather that
 of \cite{melnick} and mostly \cite{vincent}.  Those two papers present a general integrability result for Cartan geometries that will be the 
 key ingredient  
  in most of our proofs.  We summarize below the results of \cite{vincent}  (first obtained in the analytical setting in \cite{melnick}, with a 
 different approach).

Let us begin with $(M,g)$, a pseudo-Riemannian manifold of type $(p,q)$.  Let $\pi: \hm \to M$ denote the  bundle of orthonormal frames on $\hm$.  This 
 is a principal $\OO(p,q)$-bundle over $M$, and it is classical (see \cite{kobanomi}[Chap. IV.2 ])  that  the Levi-Civita connection associated to
  $g$ can be interpreted as an Ehresmann connection $\alpha$ on $\hm$, with values in the Lie algebra $\oo(p,q)$.  Let $\theta$ be the soldering 
  form on
   $\hm$, namely the $\RR^n$-valued $1$-form on $\hm$, which to every $\xi \in T_{\hx}\hm$  associates the coordinates of the vector 
   $\pi_*(\xi) \in T_xM$
    in the frame $\hx$.  The sum $\alpha + \theta$ is a $1$-form $\omega: T\hm \to \oo(p,q) \ltimes \RR^n$ called the {\it canonical Cartan 
    connection} associated 
    to $(M,g)$.  
    
    Pseudo-Riemannian structures of type $(p,q)$ are thus {\it Cartan geometries} modelled on the flat, type-$(p,q)$  Minkowski space
     $\RR^{p,q}=\OO(p,q)\ltimes \RR^n/\OO(p,q)$.  The setting of Cartan geometries being really 
     convenient for the kind of problems we are interested in, we give here the general definition.
     
     A Cartan geometry $(M,{\mathcal C})$ modelled on a homogeneous space $X=G/P$ is the data of a triple $(M,\hm,\omega)$ where $M$ is a manifold, 
     $\pi: \hm \to M$ is a $P$-principal bundle over $M$ and $\omega$, the {\it Cartan connection}, is a $1$-form on $\hm$ 
     with values in the Lie algebra $\lieg$.  There are moreover extra properties satisfied by $\omega$. 
     
     - First, for every $\hx \in \hm$,
      $\omega_{\hx}: T_{\hx}\hm \to \lieg$ is an isomorphism of vector spaces.
      
      - Moreover, the form $\omega$ is $P$-equivariant (where $P$ acts on 
      $\lieg$ via the adjoint action). 
      
      Beside pseudo-Riemannian metrics (which, as we just saw,  correspond to 
      $G=\OO(p,q) \ltimes \RR^n$, $n=p+q$, and $P=\OO(p,q)$ when the type of the metric is $(p,q)$), quite a lot of other interesting 
      geometric structures (linear connections, projective structures, conformal structures
       of dimension $\geq 3$ etc....) fit into this framework (see \cite[Chap. 4]{cap} for an extensive discussion of examples). 
       The reader wanting to learn 
       more about Cartan geometries will find modern and very comprehensive 
        introductions in \cite{cap} or \cite{sharpe}.
  
 \subsection{Generalized curvature map}
 \label{sec.curvaturemap}
 We assume now that the structures considered are of class $C^{\infty}$. The {\it curvature} of the Cartan connection $\omega$ is a $2$-form $K$ on $\hm$, with values in $\lieg$.  If $X$ and $Y$ are 
 two vector fields on $\hm$, it is given by the relation:
 $$ K(X,Y)=d \omega(X,Y) +[\omega(X),\omega(Y)].$$
 Because at each point $\hx$ of $\hm$, the Cartan connection $\omega$ establishes an isomorphism between $T_{\hx}\hm$ and $\lieg$, 
 it follows that
  any $k$-differential form on $\hm$, with values in some vector space ${\mathcal W}$,  can be seen as a map from $\hm$ to 
  $ \Hom( \otimes^k  \lieg,{\mathcal W})$.  This remark applies in particular for the curvature form, and we get  a curvature map 
  $\kappa: \hm \to {\mathcal W}_0$, where the vector space ${\mathcal W}_0$ is $\Hom(\wedge^2(\lieg/\liep);\lieg)$ (the curvature
   is antisymmetric and vanishes when one argument is tangent to the fibers of $\hm$).  
   
  We can now differentiate $\kappa$, getting a map $D \kappa : T\hm \to {\mathcal W}_0$.  Our previous remark allows to  
    see $D\kappa$  as a map  $D \kappa : \hm \to {\mathcal W}_1$, with ${\mathcal W}_1=\Hom(\lieg,{\mathcal W}_0)$. 
    Applying this procedure $r$ times, we define 
    inductively the $r$-derivative of the curvature $D^r \kappa : \hm \to \Hom( \otimes^r  \lieg,{\mathcal W}_r)$ (with ${\mathcal W}_r$ defined inductively
     by ${\mathcal W}_r=\Hom(\lieg, {\mathcal W}_{r-1})$).  
    
    Let us now set $m = \dim G$.  The {\it generalized curvature map} of the Cartan geometry $(M, {\mathcal C})$
    is the map  $\dk =(D\kappa,\ldots,D^{m+1}\kappa)$.  The $P$-module ${\mathcal W}_{m+1}$ will be rather denoted ${\mathcal W}_{\dk}$ in the following.
    
    \begin{remarque}
    Note that in the case of pseudo-Riemannian structures, the generalized curvature map encodes the  first $m+1$ covariant derivatives
     of the Riemann curvature tensor $R$.
    \end{remarque}
    
    \subsection{Integrating finite order Killing vectors}
    \label{sec.killinggenerators}

    There is a natural notion of  local isometry of a Cartan geometry $(M,{\mathcal C})$, as a local diffeomorphism $f:U \to V$ between open sets $U$
     and $V$ of $M$, which can be lifted to a local diffeomorphism of $\hat{M}$ satisfying $f^* \omega = \omega$.  In a same way, a local Killing field
      on some open subset $U \subset M$ is a vector field which can be lifted to a local vector field of $\hm$ satisfying $L_X \omega=0$. Obviously,
       those notions coincide with the classical notion of local isometry, and local Killing field when our Cartan geometry is
        defined by a pseudo-Riemannian metric on a manifold.
      
    For each integer $r \geq 1$, one defines $\Ker(D^r\kappa(\hx))$ as the vector subspace of $\lieg$ comprising
     all $\xi \in \lieg$ such that $D^r\kappa(\hx)(\xi)=0$  (recall that $D^r\kappa(\hx)$ is a linear map from $\lieg$ to ${\mathcal W}_{r-1}$).  
      
      It is clear that if $X$ is a local Killing field on $\hm$ (namely $X$ satisfies $L_X \omega=0$), 
      then all the maps $D^r\kappa$ are constant along  the orbits
       the local flow $\varphi_X^t$.  Thus $\omega(X(\hx)) \in \Ker(D^j \kappa(\hx))$ for every $j \geq 1$, 
       and every $\hx$ where $X$ is defined.  This leads naturally to a punctual notion of a {\it Killing generator of order $r$}
       at $\hx$, as a 
        vector $\xi \in \lieg$, such that $\xi$ belongs to $\bigcap_{j=1}^rKer(D^j\kappa(\hx))$.  We call $\operatorname{Kill}^r(\hx)$ the vector subspace
         of Killing generators of order $r$ at $\hx$.  
          If $m= \dim G$,  we will  note $\operatorname{Kill}^{\dk}(\hx)$ instead of  $\operatorname{Kill}^{m+1}(\hx)$. 
 
 A natural question  is now: {\it When is a Killing generator of order $r$ at $\hx$ the evaluation of an actual local 
 Killing field around $\hx$?}

 \subsubsection{The integrability theorem}   
 \label{sec.integtheorem}
 Motivated by the previous question, one defines {\it the integrability locus of $\hm$}, denoted $\hm^{\rm int}$, as follows.  A point 
  $\hx \in \hm$ belongs to 
 $\hm^{\rm int}$ if  for every 
 $\xi \in \operatorname{Kill}^{\dk}(\hx)$, there exists a Killing field $X$ defined in a neighborhood of $\hx$, and such that $\omega(X(\hx))=\xi$.  It is easily
  checked that $\hm^{\rm int}$ is a $P$-invariant set, and {\it the integrability locus of $M$}, denoted $\mint$, is just the projection of 
  $\hm^{\rm int}$ on $M$.
  Since the dimension of $\operatorname{Kill}^{\dk}(x)$ can only decrease locally, and because for every Killing field $X$, 
  the vector $\omega(X(\hx))$ belongs to  $ \operatorname{Kill}^{\dk}(\hx)$, one gets that $\mint$ is an open subset of $M$.

 At first glance, the integrability locus $\hm$ might be empty.  It turns out that it is actually dense. It seems that the first result of this kind appears
  in \cite[Th. 12]{nomizu}, for Killing generators of infinite order (namely belonging to $\operatorname{Kill}^r$ for all $r \geq 1$).  
   This  result was greatly generalized by M. Gromov in \cite{gromov}.  One crucial improvement in Gromov's approach is that it is enough to 
   consider only  Killing generators of {\it finite} order 
   (order which is moreover  controlled by the dimension of $M$ and by the nature of the geometric structure).
   The precise statement we will need in this
  article is the following (compare to \cite[Corollary 1.6.C]{gromov}, \cite[Theorem 3.11]{melnick}, \cite[Theorem 2]{vincent}):

   \begin{theoreme}[Integrability theorem]
    \label{thm.integrabilite}
    Let $(M,{\mathcal C})$ be a smooth Cartan geometry.
    The integrability locus $\mint$ coincides with  the subset of $M$ where the rank of the map $\dk$ is locally constant. 
    In particular, 
     $\mint$ is a dense open subset of $M$.
     \end{theoreme}

  Observe  that because $\dk$ is a $P$-equivariant map, the rank of $\dk$  is constant along the fibers of $\pi : \hm \to M$.
   Hence, it makes sense to speak about the rank of $\dk$ at a point $x \in M$.  More generally, we will allow in the following the notation
    $\dk(x)$ for $x \in M$, meaning by this the $P$-orbit in ${\mathcal W}_{\dk}$ of $\dk(\hx)$, for $\hx$ any point in the fiber of $x$.
  
  Although it can be easily derived  from the proofs given in \cite{vincent}, Theorem \ref{thm.integrabilite} does not appear with 
  this precise statement. For the sake of completeness, we will explain at the end of this article (Section \ref{sec.annexeA}) how to deduce
   Theorem \ref{thm.integrabilite} from \cite{vincent}. 
   
  \subsubsection{Components of the integrability locus and Kill$^{\rm loc}$-algebra} 
  \label{sec.intlocus}
  
  The dimension of the Lie algebra $\kiloc(x)$  can not  decrease locally, while the dimension of the vector space of $\operatorname{Kill}^{\dk}(x)$ (the corank of 
  $\dk$ at $x$) can not 
  increase locally.  It follows from Theorem \ref{thm.integrabilite} that on $\mint$, the Lie algebra $\kiloc(x)$ is locally constant. Hence in restriction
   to $\mint$, $\kiloc(x)$ behaves as if the structure was analytic.  More precisely
   $\mint$ splits into a union of connected components $\bigcup {\mathcal M}_i$, and on each ${\mathcal M}_i$ the Lie algebra $\kiloc(x)$ is the same
    for all $x \in {\mathcal M}_i$. We will  denote it by $\kiloc({\mathcal M}_i)$ in what follows.

%
     
    Here is another straithforward consequence of Theorem \ref{thm.integrabilite}. Assume that $x \in \mint$, 
    $\hx$ is in the fiber of $x$, and that $\{ e^{t \xi} \}_{t \in \RR}$ is a
      $1$-parameter group of $P$ fixing $\dk(\hx)$.  We get immediately that $\xi$ belongs to $\operatorname{Kill}^{\dk}(\hx)$, so that
        by Theorem \ref{thm.integrabilite}, $\omega^{-1}(\xi)$ is the evaluation at $\hx$ of a local Killing field.  This Killing field being vertical
         at $\hx$, it corresponds to a local Killing field  vanishing at $x$.  We denote in the following by ${\mathfrak I}_x$  the isotropy algebra at $x$, namely
          the Lie algebra of local Killing fields defined in a neighborhood of $x$ and vanishing at $x$.
  \begin{fact}
   \label{rk.isotropie}
   If $x \in \mint$, then the isotropy algebra ${\mathfrak I}_x$ is isomorphic to the Lie algebra of the stabilizer of $\dk(x)$ in $P$.
  \end{fact}

  Theorem  \ref{thm.integrabilite} ensures that if $x \in \mint$, $\dim \kiloc(x)$ coincides with the corank of $\dk$ at $x$. Let us call
  ${\mathcal O}_x^{\dk}$ the $P$-orbit of $\dk(\hx)$  (it does not depend on the choice of $\hx$ in the fiber of $x$, hence the notation). By Fact 
  \ref{rk.isotropie}, $\dim {\mathcal O}_{x}^{\dk} = \dim P - \dim \mathfrak{I}_x$.  We conclude that for every $x \in \mint$ :
    
  \begin{equation}
   \dim \kiloc(x) - \dim \mathfrak{I}_x=\dim M+ \dim {\mathcal O}_{x}^{\dk} - \rk(\dk)(x). \nonumber
  \end{equation}

%
%
%
%
   

\section{From quasihomogeneity to local homogeneity}
\label{sec.criteres}
 Using Theorem \ref{thm.integrabilite}, we are going to isolate some general situations where quasihomogeneity implies 
 local homogeneity.  We give below two criteria, which will be further implemented to prove Theorem \ref{thm.quasihomogene}.

\subsection{A first criterion for general Cartan geometries}
The first criterion we give is  very easy, and applies to any Cartan geometry.  That's why  we state it with this generality. We recall from the previous
 section that any Cartan geometry $(M,{\mathcal C})$ modelled on some space $X=G/P$ admits a $P$-equivariant generalized curvature map 
 $\dk : \hm \to {\mathcal W}_{\dk}$, where ${\mathcal W}_{\dk}$ is some vector space bearing a linear $P$-action. 
 
 We use in the following the notations introduced in   Section 
  \ref{sec.integration}.

\begin{proposition}
 \label{prop.critere1}
 Let $(M,{\mathcal C})$ be a connected Cartan geometry of algebraic type modelled on $X=G/P$, and $\dk$ its generalized curvature map.  
 Let $x_0 \in M$ be a point satisfying the three following properties:
 \begin{enumerate}
  \item{The Kill$^{\rm loc}$-orbit of $x_0$ is open (in particular $x_0 \in \mint$).}
  \item{The dimension of $\kiloc(x_0)$ is minimal among all the dimensions of $\kiloc(y)$, when $y$ ranges over $\mint$.}
  \item{The $P$-orbit of $\dk(x)$ is closed in ${\mathcal W}_{\dk}$.}
 \end{enumerate}
Then $(M,{\mathcal C}) $ is locally homogeneous. 
\end{proposition}

For every  $\hx \in \hm$, let us denote by ${\mathcal O}_{\hx}^{\dk}$ the orbit of  $\dk(\hx)$ in ${\mathcal W}_{\dk}$ (this orbit is the same
for all the points in a same fiber, so we will also use the notation ${\mathcal O}_{x}^{\dk}$).  Recall from Section \ref{sec.intlocus} that
 for every $x \in \mint$, we have the equality~:
  \begin{equation}
   \label{eq.dim}
   \dim \kiloc(x) - \dim \mathfrak{I}_x=\dim M+ \dim {\mathcal O}_{x}^{\dk} - \rk(\dk)(x).
  \end{equation}

The second hypothesis of the proposition says that 
$\rk(\dk)(y) \leq \rk(\dk)(x_0)$ for every $y \in \mint$. We thus will write in the following $r_{max}$ instead of $\rk(\dk)(x_0)$. 
Equation (\ref{eq.dim}), together with the first hypothesis of the proposition says that   
 $\dim {\mathcal O}_{x_0}^{\dk}=r_{max}$. 

Proposition \ref{prop.critere1} will be proved if we show that the Kill$^{\rm loc}$-orbit  ${\mathcal O}_{x_0}^{\rm loc}$, which is open by assumption, 
 is also closed. 
To do this, let us pick $x$ in the boundary of ${\mathcal O}_{x_0}^{\rm loc}$.  The $P$-orbit
  ${\mathcal O}_{x}^{\dk}$ is then in the closure of  ${\mathcal O}_{x_0}^{\dk}$, hence coincides with ${\mathcal O}_{x_0}^{\dk}$
 by the third hypothesis.  By $P$-equivariance of $\dk$, 
 the rank $\rk(\dk)(x)$ is always at least the dimension of ${\mathcal O}_{x}^{\dk}$.  This yields the inequality
 $\rk(\dk)(x) \geq r_{max}$, which must actually be an equality by definition of $r_{max}$.
  We infer that the rank of $\dk$ is locally constant around $x$, implying $x \in \mint$ by Theorem \ref{thm.integrabilite}.  
  Equality (\ref{eq.dim})
   thus holds at $x$, leading to  
   $$\dim \kiloc(x) - \dim \mathfrak{I}_x=\dim M+ \dim {\mathcal O}_{x_0}^{\dk} - r_{max}=\dim {\mathcal O}_{x_0}^{\rm loc}=\dim M.$$
   This shows that  ${\mathcal O}_x^{\rm loc}$ is open, hence meets ${\mathcal O}_{x_0}^{\rm loc}$, and 
   $x \in {\mathcal O}_{x_0}^{\rm loc}$ as desired.

\subsection{Homogeneous geodesic segments, and a  second criterion}
\label{sec.critere2}
We are now going to prove a second criterion, allowing to deduce local homogeneity of a pseudo-Riemannian structure.

Before stating the criterion, we must  recall
 the definition of a {\it homogeneous geodesic segment}. If $(M,g)$ is a pseudo-Riemannian manifold, a geodesic segment
 $\gamma : (-\delta, \delta) \to M$ is said homogeneous whenever it coincides locally with
  a piece of (local) orbit of some local Killing field. 
  For our purpose, let us make the following trivial remark:  
  if some point of a homogeneous geodesic segment is contained
   in an open Kill$^{\rm loc}$-orbit ${\mathcal O}$, then all this geodesic segment is contained in ${\mathcal O}$. 
   Thus, if some locally homogeneous component of a pseudo-Riemannian
    manifold is modelled on a space admitting a lot of homogeneous geodesic segments, one can expect to saturate the component by ``broken 
    geodesic segments'', in order
     to prove that the component is actually the whole manifold $M$. 
   
   We now make this remark into something more quantitative. For every 
   $x \in M$ and  $u \in T_xM$, denote by  $\gamma_u$ the geodesic segment
  with initial datas $\gamma_u(0)=x$ and  $\gamma_u^{\prime}(0)=u$. Let us introduce 
   the set 
   $${\mathcal A}(u)= \{ \lambda \in \R \ | \ \exists X \in \kiloc(x), \   X(x)=u \ {\rm and }\ \nabla_uX(x)=\lambda u    \},$$
   and put
   
   $$ \eta(u) = \inf \{ |\lambda| \ | \ \lambda \in {\mathcal A}(u)  \}.$$
    We adopt  the convention $\eta(u)=\infty$ when  ${\mathcal A}(u)=\emptyset$. It is  readily checked that   $\gamma_u$ is homogeneous in
    a neighborhood of $x$
   if and only if ${\mathcal A}(u)\not =\emptyset$. 
 The key fact is that if $x$ belongs to the integrability locus $\mint$, we can estimate the ``length'' of this homogeneous geodesic segment 
  inside $\gamma_u$ thanks to 
  $\eta(u)$.  More precisely:

\begin{proposition}
 \label{prop.homogene}
 Let $(M,g)$ be a pseudo-Riemannian manifold.   Let $x_0 \in \mint$, $u \in T_{x_0}M$, and $\delta>0$  such that the geodesic $\gamma_u$ 
 is defined on the interval $[-\delta,\delta]$.  Assume that the dimension of 
 $\kiloc(x_0)$ is minimal among all the dimensions of $\kiloc(y)$, when $y$ ranges over  $\mint$.
 Then, if $\tau=\min(\delta,\frac{1}{|\eta(u)|})$, the geodesic segment $\gamma_u([-\tau,\tau])$ is homogeneous.  In particular, if $\eta(u)=0$,  
 $\gamma_u([-\delta, \delta])$ is homogeneous.  
\end{proposition}

Of course, when $\eta(u)=\infty$, the proposition does not give any information.  

From proposition \ref{prop.homogene}, we can derive the following criterion of local homogeneity.
 \begin{corollaire}
  \label{coro.critere2}
  Let $(M,g)$ be a connected pseudo-Riemannian manifold, and $T\mathcal{L}M$ be the bundle of lightlike tangent vectors.  Assume that there exists 
  $x_0 \in M$ such that 
  \begin{enumerate}
  \item{The Kill$^{\rm loc}$-orbit of $x_0$ is open.}
  \item{The dimension of $\kiloc(x_0)$ is minimal among all the dimensions of $\kiloc(y)$, when $y$ ranges over $\mint$.}
  \item{For every
   compact subset $K \subset T\mathcal{L}M$, there exists $M_K>0$ such that $\eta(u) \leq M_K$ for 
   every $u \in K \cap T{\mathcal O}_{x_0}^{\rm loc}$. }
  \end{enumerate}
 Then $M$ coincides with the Kill$^{\rm loc}$-orbit of $x_0$, hence is locally homogeneous.
 \end{corollaire} 
  The proof of Proposition \ref{prop.homogene} will be done in the next section.  For the moment, let us explain how we derive Corollary 
  \ref{coro.critere2} from 
  the proposition. 
   
   We want to prove that  every $y \in M$ belongs to the Kill$^{\rm loc}$-orbit ${\mathcal O}_{x_0}^{\rm loc}$.  Let us endow $M$, with an auxiliary Riemannian metric.  The norm
   of a tangent vector $u$ for this metric will be denoted by $|u|$.  Let us consider a path $\gamma$ joining $x_0$ to $y$, and a compact set 
   $S \subset M$
    containing $\alpha$ in its interior.  The set 
    $$K=\{ v \in T_xM \ | x \in S, \ g(v,v)=0, \ |v|=1  \}$$
    is a compact subset of $T\mathcal{L}M$. 
     Let $M_K$ be the constant given by the third assumption in Corollary \ref{coro.critere2}.  It is possible to join $x_0$ and $y$  by
      a  path made of broken  lightlike  geodesic segments $\gamma_{u_i} \subset Int(S)$, with $i=0,\ldots,s$, $u_i \in K$  and $\gamma_{u_i}$ 
      defined on $[0,\epsilon_i]$, 
      $\epsilon_i < \frac{1}{M_K}$.  Now, we can apply repeatedely Proposition \ref{prop.homogene} to get that all the segments $\gamma_{u_i}$
       are homogeneous, hence $y$ belongs to ${\mathcal O}_{x_0}^{\rm loc}$. 



  
    
 \subsection{Proof of Proposition \ref{prop.homogene}}   
  \label{sec.preuveproposition}
  In what follows, we use the notations of Section \ref{sec.critere2}. 
  The set ${\mathcal A}(u)$ is an affine subspace of $\R$, hence it is either $\R$ or a  point.  One infers easily the existence of 
 $X \in \kiloc(x)$
  such that $X(x_0)=u$ and $\nabla_uX(x_0)=\eta(u)X(x_0)$.  
 We consider $\hx_0 \in \hm$ in the fiber of $x_0$, and we lift $X$ to a neighborhood $\hat U$ of $\hx_0$ in $\hm$, obtaining a field 
 (still denoted $X$) which satisfies   $L_X \omega=0$. 
   For every $\hy \in \hat U$, we put
  $\xi(\hy)=\omega(X(\hy))$.  

 \subsubsection{Exponential flow in $\hm$}  
  \label{sec.exponentialflow} 
  The Cartan connection $\omega$ on $\hm$ yields the notion of $\omega$-constant vector field on $\hm$.  If $v \in \lieg$ and $\hx \in \hm$,
  we can consider
   the local flow $\varphi_v^t$ of the $\omega$-constant vector field with value $v$, around $\hx$.  When this flow is defined up to time $1$, 
   we define
    $\exp(\hx,v)$ as $\varphi_v^1$.  The exponential map at $\hx$ : $v \mapsto \exp(\hx,v)$, is thus defined on a neighborhood of $0 \in \lieg$.
    
   Observe that the projections of curves $t \mapsto \exp(\hx,tv)$ on $M$ are geodesics for the metric $g$ when $v$ is {\it horizontal}, namely 
   $v \in \RR^n \subset \oo(p,q) \ltimes \RR^n$.  For arbitrary $v$, we get by projection a wider class of curves (for instance in the flat case, we
    get all orbits under $1$-parameter groups of $\OO(p,q) \ltimes \RR^n$).

 \subsubsection{Covariant derivative seen in $\hm$}
 \label{sec.covariant} 
  The covariant derivative $\nabla_uX(x_0)$ can be interpreted in the following way in $\hm$ (see for instance 
  \cite{kobanomi}[chap. III, $\S$1 and $\S$2]).  There is a unique 
    vector $v \in \RR^n \subset \oo(p,q)\ltimes \RR^n$  satisfying $\pi_*(\omega_{\hx_0}^{-1}(v))=u$. Now, the horizontal component of 
    $(v.\xi)(\hx_0)$ (namely the component on $\RR^n \subset \oo(p,q) \ltimes \RR^n$) is again identified by $\pi_* \circ \omega_{\hx_0}^{-1}$ with a 
   vector of 
   $T_{x_0}M$ which  is precisely  $\nabla_uX(x_0)$.   
   
   Let us now decompose $\xi=\xi(\hx_0)$ into a sum $v+A$, where $v \in \R^n$ and $A \in \oo(p,q)$, and let us 
   write $\xi(t)=\omega(\xi(\exp(\hx,tv)))$.
   
    Recall the following 
   formula for the curvature
    $K$ of $\omega$
    $$ d \omega(X,V)+[\omega(X),\omega(V)]=K(X,V)$$
    where $V$ is any vector field on $\hat U$.  Applying this formula when $V$ is the $\omega$-constant vector field such that $\omega(V)=v$,
     and using the fact that the Killing vector field $X$ commutes with $V$ (because $L_X\omega=0$), one derives easily  the following 
     differential equation satisfied by $\xi(t)$:
    
    \begin{equation}
     \label{eq.killing}
     \xi^{\prime}(t)=[\xi(t),v]-\kappa(\hgamma(t))(\xi(t),v)
    \end{equation}

    Because the curvature $K$ takes its values in $\oo(p,q)$  (this is the fact that the Levi-Civita connection is torsion-free), the horizontal
     component of $\xi^{\prime}(0)$ is $[A,v]$.  Thus our hypothesis $\nabla_uX(x_0)= \eta(u) X(x_0)$  is equivalent to $[A,v]=\eta(u) v$.
  
  \subsubsection{Saturation by the exponential flow}
  \label{sec.saturation}
       We are now going to prove an important invariance property of the component of $x_0$ with respect to the exponential flow.
  
  \begin{lemme}
  \label{lem.dkconstant}
  Let us  consider $I=(\alpha,\beta) $ a maximal interval of definition of $\hgamma: t \mapsto \exp(\hx_0,t \xi)$). Then
    the differential of the map $\dk$ is constant on $\hgamma(I)$.  In particular, the projection of $\hgamma(I)$ on $M$ is homogeneous.
  \end{lemme}

  \begin{preuve}
  We denote by $D(\dk)$ the differential of the map $\dk : \hm \to {\mathcal W}_{\dk}$. As already observed in section \ref{sec.integration}, it
  can be seen as a map from 
  $\hm$ to  $ \Hom(\lieg,{\mathcal W}_{\dk})$.
   Let us consider
   $\Lambda = \{ t \in I \ | \ D(\dk)(\hgamma(t))=D(\dk)(\hx_0)  \}$.  This is clearly a closed set of $I$.
    It is also nonempty because for $t$ in a small interval around $0$, $\exp(\hx,t\xi)$ is the orbit of a Killing field, so that 
    $D(\dk)(\hgamma(t))$ is constant 
     on this interval. Now if $t_0 \in \Lambda$, $D(\dk)(\hgamma(t_0))=D(\dk)(\hx_0)$, and we point out that the dimension of 
     Recall that for any $\hx \in \hm$, $Ker(D(\dk)(\hx))$ denotes the set of vectors $v \in \lieg$ satisfying $D(\dk)(\hx)(u)=0_{\mathcal{W}}$. 
      Observe also that the dimension of $Ker(D(\dk)(\hx))$ coincides with the corank
      of $\dk$ at $\hx$.  We thus obtain that the rank of $\dk$ is the same at $\hgamma(t_0)$ and at $\hx_0$.  The second hypothesis 
      of Proposition \ref{prop.homogene} says that
     this rank is  maximal, hence locally constant.  Since $\xi$ is in the Kernel of $D(\dk)(\hx_0)$, hence in the kernel of $D(\dk)(\hgamma(y_0))$, 
      we can apply 
        Theorem \ref{thm.integrabilite}: there exists a local Killing field $Y$ around  $\hgamma(t_0)$
       such that $\omega(Y(\hgamma(t_0)))=\xi$.  It follows that $D(\dk)$ is constant on some $\hgamma([t_0-\epsilon,t_0+\epsilon))$ for $\epsilon >0$. 
       The set $\Lambda$
        is thus open, and we are done.  
        
        The arguments above  show  that for every $t_0 \in I$, there is a small segment $\hgamma(]t_0-\epsilon,t_0+\epsilon[)$ which is 
        the orbit of
    $\hgamma(t_0)$ under the local
    flow of a local Killing field $Y$.   Projecting on $M$, we get that the segment
     $\pi(\hgamma(I))$ is  homogeneous.
  \end{preuve}
  
  \subsubsection{End of the proof of Proposition \ref{prop.homogene}}
   \label{sec.finpreuve}  
     We now make the link between the projection of $\hgamma(I)$ on $M$ and  the segment $\gamma_u$  thanks to  the following
     reparametrization lemma.
  
    \begin{lemme}[Reparametrization lemma]
  \label{lem.repara}
  Let $\xi \in \lieg$.  Assume that $\xi=A+v$, with $A \in \oo(p,q)$ and $v \in \R^n$.  
  Assume moreover that $[A,v]=\eta v$.  
  Then 
  the following equality holds for every $t$ such that at least one member makes sense:

\begin{enumerate}
\item{If $\eta=0$, $\exp(\hx,t \xi)=\exp(\hx,tv).e^{tA}\,.$}
\item{If $\eta \not = 0$, $\exp(\hx,t \xi)=\exp(\hx,\frac{1}{\eta}(e^{\eta t}-1)v).e^{tA}\,.$}
 \end{enumerate}
  
  \end{lemme}
  
  \begin{preuve}
   Let us make a first observation : two curves $t \mapsto \alpha(t)$ and $t \mapsto \beta(t)$ in $\hm$ satisfy an identity of the form
    $\alpha(t)=\beta(t).p(t)$, for some curve $t \mapsto p(t)$ with values in ${\operatorname O}(p,q)$ if and only if 
    $$ \omega(\alpha^{\prime}(t))=Ad(p(t))^{-1}\,\omega(\beta^{\prime}(t))+\omega_G(p^{\prime}(t)),$$
    where $\omega_G$ stands for the Maurer-Cartan form on the Lie group $G={\operatorname O}(p,q) \ltimes \R^n$.  
    As a consequence, we have a relation of the form
    $$ \exp(\hx,t \xi)=\exp(\hx,f(t)v).p(t)$$
    for some function $f$, if and only if 
    $$ \xi=f^{\prime}(t)Ad(p(t))^{-1}.v+\omega_G(p^{\prime}(t)),$$
     which in turns is equivalent to the identity
     $$ e^{t \xi}=e^{f(t)v}p(t)$$
     in the Lie group $G$.  
     It is thus enough to check the reparametrization formulas in $G$.

   For $m \geq 1$, we have
  $$\left( \begin{array}{cc} tA & tv\\ 0 & 0 \end{array} \right )^m= 
  \left( \begin{array}{cc} t^mA^m & \, \, \, t^m \eta^{m-1}v\\ 0 & 0 \end{array} \right ).$$
   Hence
   $$e^{t \xi}= \left( \begin{array}{cc} e^{tA} & (t+\frac{t^2 \eta}{{2! }}+\frac{t^3\eta^2}{{3!}}+ \ldots)v\\ 0 & 1 
   \end{array} \right ) $$

   This shows that $e^{t \xi}=e^{tv}.e^{tA}$ if $\eta=0$, and $e^{t \xi}=e^{\frac{1}{\eta}(e^{\eta t}-1)v}.e^{tA}$ if 
    $\eta \not = 0$, as desired.
   \end{preuve}

  We now get easily Proposition \ref{prop.homogene}, observing that the geodesic $\gamma_u$ is the projection on $M$ of the curve $t \mapsto \exp(\hx_0,tv)$.
   If we are in the first case of the reparametrization lemma (namely the case $\eta(u)=0$), then we get directly that
    $\pi(\hgamma(I))$ contains $\gamma_u([-\tau,\tau])$ (here $\tau=\delta$). By the previous remarks, $\gamma_u([-\tau,\tau])$ is then homogeneous.  
    If we are the second case of Lemma \ref{lem.repara}, it is better to rewrite the reparametrization relation 
    in the following way :
    $$ \exp(\hx,sv)=\exp(\hx, \frac{1}{\eta(u)}\ln(\eta(u) s +1))e^{-\frac{1}{\eta(u)}\ln(\eta(u) s +1)A},$$
    expression which makes sense as soon as $|s| \leq \tau$ where $\tau=\min\{ \delta; \frac{1}{|\eta(u)|} \}$. 
    Again, $\pi(\hgamma(I))$ contains $\gamma_u([-\tau,\tau])$ and we are done. 
    



  
%
%
%
%

\section{Proof of Theorem \ref{thm.quasihomogene}}
\label{sec.proofquasi}
In this section, we are under the asumptions of Theorem \ref{thm.quasihomogene}.  The manifold $(M,g)$ is $3$-dimensional and Lorentzian, 
and on a dense
 open subset, the dimension of the Lie algebra of local Killing fields is at least $4$-dimensional.

 \subsection{Hyperbolic, elliptic and parabolic components}
 
 By the general discussion made in Section \ref{sec.intlocus}, the integrability locus $\mint$  is a disjoint union of connected components: 
 $\mint=\bigcup {\mathcal M}_i$.  On each component ${\mathcal M}$, there is a well-defined Lie algebra $\kiloc({\mathcal M})$ 
 of local Killing fields, which by our assumption 
 is at least
  $4$ dimensional.  Hence, for every $x \in {\mathcal M}$, the stabilizers of the ${\operatorname O}(1,2)$-orbit $\dk(x)$ are always $1$ or $3$-dimensional, the points where
   they are $3$-dimensional being those at which the sectional curvature is constant. Because in a finite dimensional linear representation of 
   ${\operatorname O}(1,2)$, no point can have a stabilizer of dimension exactly $2$, and because of Fact \ref{rk.isotropie},  the possible dimension for  
   $\kiloc({\mathcal M})$
    is $4$ or $6$.  In the latter case, the sectional curvature is constant on ${\mathcal M}$.
    
    Assume now that the component ${\mathcal M}$ does not have constant sectional curvature. The dimension of $\kiloc({\mathcal M})$ is then 
    $4$, and we can still formulate the
    \begin{fact}[\cite{sorin.karin}]
    \label{fact.quasihomo}
    The component ${\mathcal M}$ is locally homogeneous.
    \end{fact}
 The proof of \cite[Lemma 4]{sorin.karin}  can be applied directly here.  Either all $x \in {\mathcal M}$ have $1$-dimensional isotropy algebra and we
  are done.  Or some point $x \in {\mathcal M}$ has ${\mathfrak I}_x$ of dimension $3$. The isotropy representation is then the $3$-dimensional irreducible
   representation of $\mathfrak{so}(1,2)$, and the $1$-dimensional Kill$^{\rm loc}$-orbit of $x$ yields  an invariant line for
   this representation : contradiction.  

 If $x \in {\mathcal M}$, the local flow defined by a generator of ${\mathfrak I}_x$  determines a $1$-parameter flow 
    $\varphi^t$ of 
     $\operatorname{SO}(T_xM)$, hence of $\SO(1,2)$. Up to conjugacy, there are three kinds of flows in $\SO(1,2)$:   hyperbolic, elliptic or 
     parabolic (unipotent).  We call the point $x$ hyperbolic, elliptic or parabolic, according to the type of $\varphi^t$.
     The  component ${\mathcal M}$ itself will be called hyperbolic (resp. elliptic, resp. parabolic) if it 
      contains an hyperbolic (resp. elliptic, resp. parabolic) point.  
     This terminology might seem inadequate, since for the moment, nothing prevents ${\mathcal M}$ to be hyperbolic, elliptic and parabolic
     at the same time.
       We will see in the sequel that this actually does not happen.

%

\subsection{Hyperbolic or elliptic components imply local homogeneity}
 \label{sec.hyperbolicelliptic}
 
 Let us prove Theorem \ref{thm.quasihomogene} when there exists a  component ${\mathcal M} \subset \mint$ which is hyperbolic (resp. elliptic). Let us
  pick $x_0 \in {\mathcal M}$ which is  hyperbolic or elliptic.  We are now in position to apply our first homogeneity criterion, given by
   Proposition \ref{prop.critere1}.  The two first assumptions of  Proposition \ref{prop.critere1} are fulfilled by $x_0$ (under the hypothesis 
   of Theorem \ref{thm.quasihomogene}).  To check that the third 
    condition  is also satisfied, we remark that in a finite dimensional representation of ${\operatorname O}(1,2)$, an orbit is closed if and only if
     the corresponding orbit of $\SO^o(1,2)$ is closed.  Now, $\SO^o(1,2)$ is isomorphic to $\PSL(2,\RR)$, so that the third condition follows from the 
      next 
    lemma.

 \begin{lemme}
\label{lem.orbitehyper}
Let $\rho: \PSL(2,\RR) \to \GL(V)$ be a finite dimensional representation of $\PSL(2,\RR)$.  Let $v \in V$ be a vector the stabilizer of
 which is an hyperbolic or elliptic   $1$-parameter group of $ \PSL(2,\RR)$.  Then the $\PSL(2,\RR)$-orbit of $v$  
 is closed in  $V$.
\end{lemme}

This fact is certainly standard.  For the sake of completeness, we will provide a proof at the end of the paper (see Section \ref{annexeB}).

\subsection{The geometry of parabolic components}

The previous section reduces the proof of Theorem \ref{thm.quasihomogene} to the case  when all the components of the integrability locus either
 have constant sectional  curvature, or are parabolic.    The general strategy to prove Theorem \ref{thm.quasihomogene}  will be to understand
   in more details the homogeneous spaces on which parabolic components are modelled, then to observe that those components have 
   quite a lot  of homogeneous geodesic segments,  and finally to show that Corollary  \ref{coro.critere2} can 
   be applied.
 
The main result we have in mind in this section is the description of the possible Kill$^{\rm loc}$ algebras for a parabolic component.


\begin{proposition}
  \label{prop.type.algebrique}
  Let $\mathcalm$ be a component of the integrability locus $\mint$ of a $3$-dimensional Lorentz manifold $(M,g)$.  
  We assume that $\mathcalm$ is 
  parabolic.
   Then:
  \begin{enumerate}
   \item{If the scalar curvature of $\mathcalm$ is $0$, the Lie algebra $\kiloc(\mathcalm)$ is isomorphic to 
   $\RR \ltimes_h \heis,$ $\RR \ltimes_e \heis$ or 
   $\RR \ltimes_{\nu} \heis$, with $\nu \in \RR^*$.}
   \item{If the scalar curvature is nonzero on $\mathcalm$, then  $\kiloc(\mathcalm)$ is isomorphic to $\sld \oplus \RR$.}
  \end{enumerate}
\end{proposition}
 
 The notations used in this proposition will be explained in Section \ref{sec.heisenberg} below.  Proposition \ref{prop.type.algebrique} yields a partial
  description of $4$-dimensional Lie algebras admitting a local action on a $3$-dimensional Lorentz manifolds, which is a question of independent 
  interest.  The full classification can be
   carried out with the same methods, eventhough we won't present it here.

\subsubsection{Some extensions of $\heis$ by $\RR$}
\label{sec.heisenberg}
Let us make some algebraic preliminaries about $4$-dimensional Lie algebras which are a semidirect product
$\RR \ltimes \heis$, where $\heis$ stands for the $3$-dimensional Heisenberg Lie algebra.

Let us consider  a Lie algebra $\lieg \simeq \RR \ltimes \heis$
 with basis $X,Y,Z,T$, where $X,Y,Z$ generates the Lie subalgebra $\heis$, with relation $[X,Y]=Z$ (and all  other
  brackets involving $X,Y,Z$ zero), and $\ad(T)$ is given in the basis $X,Y,Z$ by a matrix
  $$ \left( \begin{array}{cc}  A & 0\\ 0 & \operatorname{Tr}(A) \end{array}  \right).$$
  If we replace  $T$ by $\alpha T$, $\alpha \not =0$, we won't affect the isomorphism type of our Lie algebra.  The same is true if  for some 
  $P \in \GL(2,\RR)$, we replace $X,Y$ by $P.X,P.Y$, and $Z$ by $\operatorname{det}(P).Z$.  Thus, we can replace $A$ by any multiple $\alpha A$, 
  $\alpha \not =0$, or a conjugate $PAP^{-1}$ with $P \in \GL(2,\RR)$.  We end up with the following few cases:
  
  \begin{enumerate}
  \item{{\it Unimodular case}. If the trace of $A$ is $0$,   we may choose 
  $$A=0,   \ A=\left( \begin{array}{cc}  1& 0\\ 0 & -1 \end{array}  \right), \  
   A=\left( \begin{array}{cc}  0& 1\\ 0 & 0 \end{array}  \right) \  {\rm or} \  A=\left( \begin{array}{cc}  0& 1\\ -1 & 0 \end{array}  \right).$$
   The corresponding Lie algebras are denoted respectively :\\
   $\RR \times \heis$, $\RR \ltimes_h \heis$, $\RR \ltimes_p \heis$ and $\RR \ltimes_e \heis$.}
  \item{{\it Scalar case}. If $A$ is a nonzero scalar matrix, we may choose $A=\left( \begin{array}{cc}  1& 0\\ 0 & 1 \end{array}  \right)$, and we
   denote the corresponding Lie algebra by $\RR \ltimes_{s} \heis$.}
  \item{{\it Other nonunimodular cases}.  In all other cases, the matrix $A$ is not scalar and has a nonzero trace, that we might assume to be $2$. This leads to 
  three subcases
  \begin{enumerate}
   \item{Hyperbolic case: $A=\left( \begin{array}{cc}  1+\lambda& 0\\ 0 & 1-\lambda \end{array}  \right), \ \lambda>0.$}
   \item{Unipotent case: $A=\left( \begin{array}{cc}  1& 1\\ 0 & 1 \end{array}  \right)$.}
   \item{Elliptic case: $A=\left( \begin{array}{cc} \lambda & 1\\ - 1 &  \lambda \end{array}  \right), 
   \ \lambda>0$.}
  \end{enumerate}
}
  \end{enumerate}

Observe that in case $(3)$ above (namely $A$  is not scalar and $\operatorname{Tr}(A) \not =0$),  
the number
$\nu=\frac{\operatorname{det}(A)}{\operatorname{Tr}(A)^2}$ is a complete  invariant of isomorphism classes. We will then denote by
$\RR \ltimes_{\nu} \heis$ this last family of algebras. The parameter $\nu$ can take any real value. For $\nu<1/4$, 
we are in the hyperbolic case above, the parabolic and elliptic cases corresponding to $\nu=1/4$ and $\nu >1/4$.

%
%
%

\subsubsection{The curvature module}

We consider on  $\RR^3$ the Lorentzian form, with matrix in a   basis $e,h,f$ given by:
$$ J=\left(  \begin{array}{ccc}    0&0&1\\ 0&1&0 \\ 1&0&0\\ \end{array}  \right).$$
We call ${\operatorname O}(1,2)$ the subgroup of  $\GL(3,\RR)$ preserving the bilinear form determined by $J$. 
Its Lie algebra is denoted by   $\oo(1,2)$, and admits the following basis :
$$ E=\left(  \begin{array}{ccc}   0&1&0\\ 0&0& -1\\ 0&0&0\\ \end{array}  \right),  
H= \left(  \begin{array}{ccc} 1&0&0\\ 0&0&0 \\ 0& 0& -1   \end{array}  \right)$$
$$ F= \left(  \begin{array}{ccc}   0&0&0\\ 1&0&0 \\ 0&-1& 0 \end{array}  \right).$$
We thus have the comutation relations $[H,E]=E, [H,F]=-F$ and $[E,F]=H$.

Let $(M^3,g)$ be  $3$-dimensional Lorentz manifold, and denote by  $\hm$ its bundle of orthonormal frames.  At each $\hx \in \hm$, the 
 curvature $\kappa(\hx)$ is an element of $\Hom(\wedge^2(\RR^3),\oo(1,2))$. Because of Bianchi's identities,
 the curvature module is actually a $6$-dimensional
  submodule of $\Hom(\wedge^2(\RR^3),\oo(1,2))$.
Chosing $e \wedge h$, $e \wedge f$, $h \wedge f$ as a basis for $\wedge^2(\RR^3)$, and $E,H,F$ as a basis for $\oo(1,2)$, 
an element of $\Hom(\wedge^2(\RR^3),\oo(1,2))$ is merely given by a $3 \times 3$ matrix.   The action of ${\operatorname O}(1,2)$ on 
$\Hom(\wedge^2(\RR^3),\oo(1,2))$ corresponds to  the conjugation on matrices.

Scalar matrices are ${\operatorname O}(1,2)$-invariant, and form a $1$-dimensional irreducible submodule 
(corresponding to constant sectional curvature).

The other irreducible submodule of the curvature module is $5$-dimensional, spanned by the matrices:

$$ \left(  \begin{array}{ccc}   0&0&1\\ 0&0&0\\ 0&0&0\\  \end{array}  \right), 
 \left(  \begin{array}{ccc}   0&1&0\\ 0&0&1\\ 0&0&0\\ \end{array}  \right), 
 \left(  \begin{array}{ccc}  1&0&0\\ 0&-2 & 0\\0& 0& 1 \\  \end{array}  \right),$$
$$ \left(  \begin{array}{ccc}  0&0&0\\ 1&0 &0 \\ 0 & 1 & 0\\  \end{array}  \right), \left(  \begin{array}{ccc}    
0&0&0 \\ 0&0&0 \\1&0&0 \\ \end{array}  \right).$$

We call $\kappa_0$ the element of $\Hom(\wedge^2(\RR^3),\oo(1,2))$  corresponding to the identity matrix, namely
$\kappa_0$ maps $e \wedge h$ to $E$, $e \wedge f$ to $H$ and  $h \wedge f$ to $F$.  We also call $\kappa_1$ 
the element of  $\Hom(\wedge^2(\RR^3),\oo(1,2))$ 
corresponding to the matrix $\left(  \begin{array}{ccc}   0&0&1\\ 0&0&0\\ 0&0&0\\  \end{array}  \right)$.  

The two dimensional vector space spanned by $\kappa_0$ and $\kappa_1$ is the set of fixed points of the action of $\{ e^{tE} \}_{t \in \RR}$ on 
the curvature module.
  


  
\subsubsection{Equation of Killing fields}
 
  Cartan's formula for the Lie derivative $L_X =\iota_X \circ d + d \circ \iota_X$ yields that whenever $U$ and $V$ are two 
 local Killing fields on an open subset of  $\hm$, the following relation holds:
 \begin{equation}
 \label{eq.killing}
 \omega([U,V])=\kappa(\omega(U) \wedge \omega(V))-[\omega(U),\omega(V)])
 \end{equation}
 
 Let ${\mathcal M}$ be a parabolic component, and $x$  a point of $\mathcalm$ having an open Kill$^{\rm loc}$-orbit (see Fact \ref{fact.quasihomo}). 
  Because  $\mathcalm$ is parabolic, there exists 
 a local Killing field $Y$, vanishing at $x$, and generating
  a parabolic flow of $\SO(1,2)$ in $T_xM$.  Thus, there exists $\hx$ in the fiber of $x$ such that after lifting $Y$ into a 
  Killing field around $\hx$, we have $\omega_{\hx}(Y(\hx))=E$.  We can also find $3$ other local Killing fields $T,X,Z$  such that at $\hx$:
 \begin{equation}
  \label{eq.omega}
  \omega(Z)=e+\mu H + \zeta F,
   \end{equation}
 \begin{center}$\omega(X)=h+\alpha H +\theta F \text{   and     } \omega(T)=f+ \beta H + \gamma F.$ \end{center}

 Observe also that at
 $\hx$, the curvature $\kappa(\hx)$ is invariant by the action of $\{e^{tE}\}$, hence is of  the form $\sigma(x) \kappa_0+b(\hx)\kappa_1$. 
  We wrote $\sigma(x)$ instead of $\sigma(\hx)$  since this number does not depend on the point $\hx$ in the fiber of $x$.  Actually, 
  $x \mapsto \sigma(x)$ coincides with (a constant multiple of) the scalar curvature.  Since ${\mathcal M}$ is locally homogeneous 
  (see Fact \ref{fact.quasihomo}), $\sigma$ is constant on ${\mathcal M}$. To simplify the notations, we will write $b$ instead of $b(\hx)$  (keeping in
   mind, though, that it depends on $\hx$).  
   This yields  
 the relations
 \begin{equation}
  \label{eq.courbure}
  \kappa(e \wedge h)= \sigma E, \ \kappa(e \wedge f)=\sigma H, \ \kappa(h \wedge f)=bE+\sigma F.
 \end{equation}

 Relation (\ref{eq.killing}) expressed for the $4$ Killing fields $X,Y,Z,T$ provides $6$ equations in $\mathfrak{o}(1,2)\ltimes \RR^3$  
  involving $\alpha,\beta,\gamma,\mu,\theta,\zeta$, as well as the two auxiliary parameters $\sigma$ and $b$.
    We are going to solve those
   equations explicitely below.
%
%
%
%
%
%

\subsubsection{Bracket relations and proof of proposition \ref{prop.type.algebrique}}
\label{sec.proofalgebric}
Let us introduce  the vector space 
    $$\mathcal{E}={\operatorname{Span}}(\omega(X(\hx)),\omega(Y(\hx)),\omega(Z(\hx)),\omega(T(\hx))).$$
Relations (\ref{eq.courbure}) yield  the following identities at $\hx$ :

$$\begin{array}{ccccccc}
 \omega([T,Y]) &=& -[\omega(T),\omega(Y)] &=& -h-\beta E + \gamma H &=& \omega(-X-\beta Y)+(\alpha+\gamma)H+\theta F,\\
 \omega([Y,X]) &=&-[\omega(Y),\omega(X)] &=& - e + \alpha E &=& \omega(\alpha Y -Z)  +\mu  H,\\
 \omega([Z,Y]) &=&  -[\omega(Z),\omega(Y)] &=&  -[e+ \mu H + \zeta F,E]&=& -\mu\omega(Y) + \zeta H.
\end{array}$$

Because all the expressions above must belong to ${\mathcal E}$, we infer  
 \begin{equation}
 \label{eq.ty}
 \alpha=-\gamma \text{ and }\theta=\zeta=\mu=0.
 \end{equation}
 Moreover, since a local Killing field $U$ satisfying $\omega(U(\hx))=0$ must be zero 
 (a local flow of   isometries must act freely on the bundle of orthonormal frames $\hm$), we also get the identities 
 \begin{equation}
  \label{eq.bracket1}
  [T,Y]=-X-\beta Y, \ [X,Y]= Z - \alpha Y  \text{ and }  [Y,Z]=0.
 \end{equation}

Relations (\ref{eq.ty}) allow to update (\ref{eq.omega}), and we get at $\hx$ the new identities :
$$ \omega(Z)=e, \omega(Y)=E \text{ and } \omega(X)=h+\alpha H .$$
We then write $3$ more equations:
 \begin{eqnarray}
 \omega([Z,X]) &=&   \ka(e \wedge h)-[e, h + \alpha H ] = \sigma E + \alpha e = \omega(\alpha Z + \sigma Y )\nonumber\\
 \omega([T,X]) &=&   \kappa(f \wedge h)+[h+\alpha H, f +  \beta H -\alpha F]\nonumber\\
   &=&  -2 \alpha f -b E + (\alpha^2-\sigma)F \nonumber\\
  &=& \omega( -2 \alpha T  -b Y)+2 \alpha \beta H -(\alpha^2+\sigma)F \nonumber\\
 \omega([T,Z])&=&\kappa(f \wedge e)+[\omega(Z),\omega(T)] = -\sigma H+[e,f+\beta H -\alpha F]\nonumber\\ 
 &=&\omega(\alpha X - \beta Z)-(\alpha^2 + \sigma)H.\nonumber
 \end{eqnarray}

 which imply, by the same arguments as above
\begin{equation}
 \label{eq.relations}
 \alpha \beta =0 \text{ and } \alpha^2 = -\sigma.
\end{equation}
as well as $3$ more bracket relations:
\begin{equation}
\label{eq.bracket2}
[X,Z]=-\alpha Z -\sigma Y, [T,X]=-2 \alpha T - b Y, \text{ and }  [T,Z]=\alpha X - \beta Z.
\end{equation}
This leads to  two cases according to the value of the scalar curvature $\sigma$.

\subsubsection*{Case of zero scalar curvature}     Because of relation (\ref{eq.relations}), the hypothesis $\sigma=0$ 
 also implies $\alpha=0$.  Relations (\ref{eq.bracket1}) and (\ref{eq.bracket2}) imply $[Y,Z]=[X,Z]=0, [X,Y]=Z$, so that the Lie algebra $\lieh$
generated by $X,Y,Z$ is  isomorphic to $\heis$, and is an ideal in $\kiloc(x)$. Moreover the adjoint action of $T$ on $\lieh$,
expressed in the basis $Y,X,Z$ is given by the matrix
$$ \ad(T)= \left (   \begin{array}{ccc} -\beta & -b & 0\\ -1 & 0 & 0\\ 0 & 0 & -\beta \\ \end{array} \right ).$$    
The Lie algebra $\kiloc(x)$ is a semi-direct product $\RR \ltimes \heis$ described by the matrix $A=\left( 
\begin{array}{cc} -\beta & -b\\ -1 & 0 \end{array}\right)$, as explained in Section \ref{sec.heisenberg}.  We see that the matrix $A$ 
is never scalar.  Moreover, observe that while $b=b(\hx)$ has no precise meaning on $M$, its vanishing has one. It means that the sectional curvature
    is constant at $x$.  In particular, because the Kill$^{\rm loc}$-orbit of $\hx$ is open, and  because ${\mathcal M}$ does not have constant sectional
     curvature (on any open set) by hypothesis, {\it we see that $b \not =0$}.  As a consequence, the matrix $A$ is invertible. 
     Hence  $\kiloc(x)$ is isomorphic to the Lie algebra $\RR \ltimes_h \heis$ (case $\beta=0$, $b<0$), 
$\RR \ltimes_e \heis$ (case $\beta=0$, $b>0$), or to the Lie algebra $\RR \ltimes_{\nu} \heis$, where $\nu=\frac{-b}{\beta^2}$ 
(here $\beta \not =0$) can take all nonzero  real values.

\subsubsection*{Case of nonzero scalar curvature} Relation (\ref{eq.relations}) yields $\sigma=-\alpha^2$, so that $\sigma$ turns 
out to be negative. Also, (\ref{eq.relations}) shows that $\beta=0$. The bracket relations become:

$$ [T,Z]=\alpha X \  \ [T,Y]=-X \  \ [T,X]=-2 \alpha T -b Y,$$
$$ [X,Z]=-\alpha Z + \alpha^2 Y \ \ [X,Y]=Z-\alpha Y \ \ [Z,Y]=0.$$

Let us put $Z^{\prime}=\frac{-1}{4 \alpha^2}(Z + \alpha Y)$, $Y^{\prime}=\frac{1}{4 \alpha^2}(Z - \alpha Y)$, 
$X^{\prime}=\frac{-1}{2\alpha}X$ and $T^{\prime}=T-\frac{b}{2}Y^{\prime}-bZ^{\prime}.$

It is easily checked that now, $Z^{\prime}$ lies in the center of $\kiloc(x)$.  Moreover:
$$  [X^{\prime},Y^{\prime}]=\frac{-1}{8 \alpha^3}[X,Z-\alpha Y]=\frac{2 \alpha}{8 \alpha^3}(Z-\alpha Y)=Y^{\prime}$$
\begin{eqnarray} 
[X^{\prime},T^{\prime}]&=&[X^{\prime},T]-\frac{b}{2}[X^{\prime},Y^{\prime}]=-T-\frac{b}{2 \alpha}Y-\frac{b}{2}Y^{\prime}\nonumber\\
&=&-T+\frac{b}{2}Y^{\prime}+b Z^{\prime}\nonumber \\
&=& -T^{\prime} \nonumber
  \end{eqnarray}

and  $[Y^{\prime},T^{\prime}]=[Y^{\prime},T]=-\frac{1}{2\alpha}X=X^{\prime}$.

Thus, the Lie algebra $\kiloc(x)$ is isomorphic to a product $\sld \oplus \RR$. 

\subsection{End of the proof of Theorem \ref{thm.quasihomogene}}

By Section \ref{sec.hyperbolicelliptic}, we can assume that all the components of $\mint$ are parabolic or of constant sectional curvature. 
If all are of constant sectional curvature, then $(M,g)$ itself has constant sectional curvature (by density of $\mint$), and $(M,g)$ is
locally homogeneous.
We thus assume that there exists   at least one  parabolic component ${\mathcal M}$, and we want to apply Corollary \ref{coro.critere2} in order to 
show that 
$M$ is locally homogeneous.  Any point $x_0 \in \mathcal{M}$ having an open orbit satisfies the   two first conditions
 of Corollary \ref{coro.critere2}. The harder part is  to check the third condition.
 
  Let us consider the closed subset $\hat{F} \subset \hm$ defined as $\hat{F}=\kappa^{-1}(\operatorname{Span}(\kappa_0,\kappa_1))$.  Because
   we excluded the presence of hyperbolic and elliptic components, we know that in each fiber of $\hm$, there is a point where $\kappa$ is 
   fixed by the $1$-parameter group $\{e^{tE}\}_{t \in \RR}$.  It follows that $\hat{F}$ projects onto $M$. By the very definition of $\hat{F}$, there exist two continuous functions $\sigma$ and $b$ on 
  $\hat F$,  such that for every $\hx \in \hat{F}$, $\kappa(\hx)=\sigma(\hx) \kappa_0 + b(\hx) \kappa_1$.  Actually, because $\sigma$ is constant 
  along the fibers, and because $M$ is locally homogeneous on a dense open set, the function $\sigma$ is locally constant on a dense open subset of $\hm$, 
  hence constant. We will  thus write $\sigma$ instead of $\sigma(\hx)$. 
  
  At each $\hx \in \hm$, there is a linear map 
  $\iota_{\hx} : \lieg \to T_xM$ defined by $\iota_{\hx}(u)=\pi_*(\omega_{\hx}^{-1}(u))$.  This map $\iota_{\hx}$ is an isomorphism between
   $\RR^3 \subset \oo(1,2)\ltimes \RR^3$ and $T_xM$.
  
  In the statement below,  we will denote by $f_t$ the vector $Ad(e^{tE}).f$.

\begin{lemme}
 \label{lem.nu}
 Let  $\mathcal M$ be a parabolic component, $x$ a point of $\mathcalm$, and $\hx$  a point of $\hat F$ projecting on $x$.  Then 
 \begin{enumerate}
  \item{If $\lieg$ is isomorphic to $\sld \oplus \RR$,  $\RR \ltimes_h \heis$ or $\RR \ltimes_e \heis$,  
   then $\eta(\iota_{\hx}(e))=0=\eta(\iota_{\hx}(f_t))$ for all $t \in \RR$.}
  \item{If $\lieg$ is isomorphic to $\RR \ltimes_{\nu} \heis$, $\nu \in \RR^*$, then $
 \eta(\iota_{\hx}(e))=0$ and for all $t \in \RR$, $\eta(\iota_{\hx}(f_t)) \leq  \sqrt{\frac{|b(\hx)|}{|\nu|}}$.}
 \end{enumerate}

 \end{lemme}
  
 \begin{preuve}
If we look at the proof of Proposition \ref{prop.type.algebrique} given in Section \ref{sec.proofalgebric},
we see that in each case (zero or nonzero scalar curvature), and for each $\hx \in \hat{F}$, we have
 three Killing fields $X$, $Y$, $T$ in $\kiloc$ satisfying:
 
 \begin{enumerate}
  \item{If $\lieg$ is isomorphic to $\sld \oplus \RR$,  $\RR \ltimes_h \heis$ or $\RR \ltimes_e \heis$,  
    $ \omega(X)=e, \omega(Y)=E, \omega(T)=f-\alpha F$, with $\alpha \in \RR$ (evaluation made at $\hx$).}
 \item{If $\lieg$ is isomorphic to $\RR \ltimes_{\nu} \heis$, $\nu \in \RR^*$, 
 $\omega(X)=e, \omega(Y)=E,  \omega(T)=f+\beta H$, with $\beta \not =0$ satisfying $\beta^2=\frac{- b(\hx)}{\nu}$ 
 (again, the evaluation is made at $\hx$). }
 \end{enumerate}
 
 This leads immediately to $\nabla_{\iota_{\hx}(e)}X(x)=0$. Because $F(f)=0$ and $H(f)=-f$, one has $\nabla_{\iota_{\hx}(f)}T(x)=0$ in the first case
  and $\nabla_{\iota_{\hx}(f)}T(x)=-\beta T(x)$ in the second one, where $\beta = \pm  \sqrt{\frac{|b(\hx)|}{\nu}}$.  Both cases 
  lead to the conclusion of the lemma for $t = 0$.
  
  Now the local flow of $Y$, denoted $Y^t$, acts on $\kiloc(x)$.  For every $t \in \RR$, $Y^t$ is defined on a small neighborhood of $x$, and
   $T_t=(Y^t)^*T$ is a local Killing field around $x$.  One checks that $T_t(x)=\iota_{\hx}(f_t)$.  Because the local flow $Y^t$ preserves 
   $\nabla$ one still has $\nabla_{\iota_{\hx}(f_t)}T_t(x)=0$ in the first case above
  and $\nabla_{\iota_{\hx}(f_t)}T_t(x)=-\beta T_t(x)$ in the second one, which establishes the lemma.  
 \end{preuve}

  
  We consider $K$ a compact subset of $T{\mathcal L}M$, projecting on a compact set $C \subset M$, and such that $C \cap M \not = \emptyset$. Let us take a compact subset $\hat C \subset \hm$
   which project on $C$.  We can saturate ${\hat C}$ by the action of $\operatorname{O}(2) \subset \operatorname{O}(1,2)$, 
  and we denote again ${\hat C}$ the compact set  obtained in this way. Observe that for every $\hx \in \hm$, the orbit of $\hx$ under 
  the right action of $\operatorname{O}(2)$ meets $\hat{F}$.  This is just because the conjugates of $\RR.E$ under 
  $\operatorname{O}(2)$ describe all the parabolic directions in $\mathfrak{so}(1,2)$ 
  ({\it i.e} directions generating a parabolic $1$-parameter group).  Thus, we can replace ${\hat C}$  by ${\hat F} \cap {\hat C}$, 
  obtaining a compact subset
   which still projects onto $C$.  
   
   Let us now endow $\RR^3$ with a norm $|.|$. By compactness of ${\hat C}$ and $K$, there exist positive $C_1$ and $C_2$ such that 
   $\iota_{\hx}^{-1}(u) \leq |C_1|$ and $|b(\hx)| \leq C_2$ for every $\hx \in \hat C$ and $u \in K$.  
   
   It is clear that  $\lim_{|t| \to \infty} |Ad(e^{tE}).f|=\infty$, so that there exists $C_3>0$ with 
   $\min\{ |e|, \inf_{t \in \R}|Ad(e^{tE}).f| \} \geq C_3$.

   Let us pick  a lightlike vector $u \in K \cap T{\mathcal M}$. Let us choose any $\hx$ in $\hat C$ that projects on $x$.
   We thus have either $\iota_{\hx}^{-1}(u)= \pm |\iota_{\hx}^{-1}(u)|e$, or
   $\iota_{\hx}(u)=\pm \frac{|\iota_{\hx}^{-1}(u)|}{|f_t|}f_t$ for some $t \in \RR$.  In the first case, Lemma \ref{lem.nu} 
   says that $\eta(u)=0$.  In the second case, we use first the fact  that $\eta$ is $1$-homogeneous to get 
   $\eta(u)=\frac{|\iota_{\hx}^{-1}(u)|}{|f_t|}\eta(\iota_{\hx}(f_t))$.  Then we  apply Lemma  \ref{lem.nu} and obtain
   $$ \eta(u) \leq \frac{ C_1\sqrt{C_2}}{\sqrt{|\nu|}C_3}.$$
   
   The third condition  of Corollary \ref{coro.critere2} is satisfied, and  the proof of Theorem \ref{thm.quasihomogene} is complete.

\section{Dense orbit implies locally homogeneous}
\label{sec.dense.orbit}

Let us now explain how one can deduce Theorem \ref{thm.orbitedense}  from  Theorem \ref{thm.quasihomogene}.  We start  with a
 smooth, closed $3$-dimensional Lorentz manifold $(M^3,g)$, and our assumption is that $\Iso(M,g)$ admits a dense orbit. Then Gromov's theorem
   \ref{thm.gromov}  ensures that the integrability locus $\mint$ is  quasihomogeneous.  Actually, because there is a dense orbit for 
   $\Iso(M,g)$, all the components of $\mint$ are pairwise isometric.  Thus the Lie algebra $\kiloc({\mathcal M})$
   does not depend on the component
    $\mathcal M$, and will be merely denoted $kiloc$.  Since we are in a quasihomogeneous situation, the dimension of this algebra is $3$, $4$ or $6$ 
    (we already noticed that $5$ is not allowed).  
      If $\kiloc$ is $6$-dimensional, then $g$ has constant sectional curvature on $\mint$, hence on $M$ by density of $\mint$  in $M$, 
      and $(M,g)$ is indeed
       locally homogeneous.  If $\kiloc$ is $4$-dimensional, then Theorem \ref{thm.quasihomogene} applies, and we conclude again that
        $(M,g)$ is locally homogeneous.  Hence, Theorem \ref{thm.orbitedense} will be proved if we can show that 
     $\lieg$ is not $3$-dimensional. This is at this point that we will use the compacity of $M$  (actually a finite volume assumption 
     would be enough).
     
     Observe that we can assume $\Iso(M,g)$ noncompact, otherwise $(M,g)$ would be directly homogeneous under
      the action of $\Iso(M,g)$. Observe also that whenever $\kiloc$ is $3$-dimensional,  the isotropy algebra $\mathfrak{I}_x$ must be trivial 
      for every $x \in \mint$. The conclusion will thus  be a direct consequence of the following general observation:
      
      \begin{proposition}
       \label{prop.stabilisateurs}
       Let $(M,g)$ be a compact pseudo-riemannian manifold. If the group $Iso(M,g)$ is noncompact, 
then for almost every $x \in \mint$, the isotropy algebra $\mathfrak{I}_x$ generates a noncompact subgroup of ${\operatorname O}(T_xM)$
      \end{proposition}

  \begin{preuve}
   We denote by $(p,q)$ the signature of $g$, $\hm$ the bundle of orthonormal frames on $M$, and  $\dk$ the generalized curvature map. 
   Recall from section \ref{sec.curvaturemap}
    that $\dk$ has range into the ${\operatorname O}(p,q)$ module ${\mathcal W}_{\dk}=\Hom( \otimes^{m+1}  \lieg,{\mathcal V})$, where 
${\mathcal V}$ is $\Hom(\wedge^2(\RR^n);\oo(p,q) \ltimes \RR^n)$, and $m=\frac{(p+q)(p+q+1)}{2}$.  
   
   Let us first recall the following recurrence theorem :

\begin{lemme}[Poincar\'e recurrence, see \cite{feres}, Theorem 2.2.6]

Let $G$ be a Lie group acting continuously on a manifold $M$, and
 preserving a finite Borel  measure $\nu$.  Then for almost every $x \in M$, there exists a sequence 
$(g_k)$ leaving every compact subset of $G$, and a sequence $(x_k)$ converging to $x$, so that $g_k.x_k$ converges to $x$.  

\end{lemme}
 
 The lemma applies to $(M,g)$ because a closed Lorentz manifold naturally defines a finite Borel measure which is invariant by isometries 
 (the density giving volume $1$ to every orthonormal frame).  The set of  points in $\mint$ which are recurrent for $\Iso(M,g)$ has 
 thus full measure in $\mint$. Let
 $x$ be such a point, and 
 $\hx \in \hn$   in the fiber of $x$.  The recurrence hypothesis means that 
    there exists $(f_k)$ 
 tending to infinity in $\Iso(M,g)$, and $(p_k)$ a sequence of ${\operatorname O}(p,q)$ such that $f_k(\hx).p_k^{-1}$ tends to $\hx$. 
 Observe that $(p_k)$ tends to infinity in ${\operatorname O}(p,q)$, because $\Iso(M,g)$ acts properly on $\hm$.  By equivariance of $\dk$, we also have
 $$p_k.\dk(\hx) \to \dk(\hx).$$
 
   The action of the algebraic group ${\operatorname O}(p,q)$ on $\mathcal{W}_{\dk}$ is linear, hence algebraic. 
   As a consequence,  all
 the orbits of ${\operatorname O}(p,q)$ are locally closed.  It is in particular the
 case of the orbit $\mathcal O$ of $\dk(\hx)$. If $I$ denotes the
  stabilizer of $\dk(\hx)$ in ${\operatorname O}(p,q)$, then the orbital map ${\operatorname O}(p,q)/I \to \mathcal O$ is an homeomorphism, when $\mathcal O$ is endowed with the topology 
  induced by that of $\mathcal{W}_{\dk}$.  As a consequence, the property $p_k.\dk(\hx) \to \dk(\hx)$ implies the existence of a sequence
   $(\epsilon_k)$ in ${\operatorname O}(p,q)$ with $\epsilon_k \to Id$ and  $\epsilon_k.p_k.\dk(\hx) = \dk(\hx)$.
  Since $(p_k)$ tends to infiny,  so does  $(\epsilon_k.p_k)$, proving that $I$ is noncompact.  Because $I$ is an algebraic group, the identity component
   $I^{o}$ is noncompact, and we conclude thanks to Fact \ref{rk.isotropie}.

  \end{preuve}

 \section{Constructing quasihomogeneous metrics of lower regularity}
 \label{sec.c1}
 Theorem \ref{thm.quasihomogene}  is stated in the framework of smooth lorentz structures.  Actually,  a closer look at the proof shows that we 
 need a regularity of the metric yielding 
  a $C^1$ generalized curvature map $\dk$. 
  Hence, because the generalized curvature map involves, for $3$-dimensional Lorentz metrics, the $6$ first covariant derivatives of the 
  curvature,
   our proof is actually available for metrics of class $C^9$.....
   
 It is likely that Theorem \ref{thm.quasihomogene} holds for metrics of lower regularity. However, we are going to exhibit $C^1$ quasihomogeneous Lorentz metrics 
  which are not locally homogeneous.  Hence, regularity at least $C^2$ is necessary to get Theorem \ref{thm.quasihomogene}.  Here is our 
  statement:

  \begin{theoreme}
  \label{theo.c1}
   There exist $C^1$  lorentzian $3$-manifolds $(M,g)$ which are quasihomogeneous but not locally homogeneous.  More precisely, for any
    pair $(\nu_1,\nu_2)$  of distinct numbers in $(- \infty, \frac{1}{4})$, one can buildt a $C^1$ Lorentz manifold $(M,g)$ having the following
     properties:
     \begin{enumerate}
      \item{There exists an open subset $U_1 \subset M$ which is locally homogeneous, with local Killing algebra isomorphic to
      $\RR \ltimes_{\nu_1} \heis$.}
      \item{There exists an open subset $U_2 \subset M$ which is locally homogeneous, with local Killing algebra isomorphic to
      $\RR \ltimes_{\nu_2} \heis$.}
      \item{The union $U_1 \cup U_2$ is dense in $M$.}
     \end{enumerate}

  \end{theoreme}

 We refer to Section \ref{sec.heisenberg} for the notation $\RR \ltimes_{\nu} \heis$, and the classification of extensions 
 of $\heis$ by $\RR$.
 
 Actually, if we drop the regularity to $C^0$, it is even possible to get quasi homogeneous Lorentz metrics on compact manifolds.
 
 \begin{theoreme}
  \label{thm.quasicompact}
  There are $C^0$ Lorentz metrics on the $3$-torus which are quasihomogeneous but not homogeneous.
 \end{theoreme}

  The constructions work as follows.  Let us endow $\RR^3$  with coordinates $(x_1,x_2,x_3)$, and let $g_0=-2dx_1dx_3+dx_2^2$  be the 
  Minkowski metric.  If $\alpha \in \RR^*$, we can define on the open set $x_3>0$
   the metric $g_{\alpha}=x_3^{\alpha}g_0$.  
   
   We introduce  the $3$ vector fields:
   $$ X=x_2 \frac{\partial }{\partial x_1} + x_3 \frac{\partial }{\partial x_2}, \ \ \ Y=\frac{\partial }{\partial x_2}, \ \ \ Z=\frac{\partial }{\partial x_1}.$$
   Those are Killing fields for the flat metric $g_0$, and because they have vanishing component along $\frac{\partial }{\partial x_3}$, they are also Killing fields
    for $g_{\alpha}$.
    
    There is a fourth vector field 
    $T=2\frac{\alpha+1}{\alpha}x_1 \frac{\partial }{\partial x_1} + x_2 \frac{\partial }{\partial x_2} 
    -\frac{2}{\alpha}x_3 \frac{\partial }{\partial x_3}$ 
     which is also Killing for $g_{\alpha}$. Indeed, the flow $\varphi^t$ generated by $T$ is the linear flow 
    $$ e^t \left(  \begin{array}{ccc}  e^{\frac{\alpha+2}{\alpha}t} & 0 & 0 \\ 0 & 1 & 0 \\ 0 & 0 & e^{-\frac{\alpha+2}{\alpha}t}
    \end{array} \right),$$
    hence $(\varphi^t)^*g_0=e^{2t}g_0$  and 
    $$(\varphi^t)^*g_{\alpha}=(e^{(1-\frac{\alpha+2}{\alpha})t}x_3)^{\alpha}e^{2t}g_0=g_{\alpha}.$$
    
    The only nontrivial bracket relations between $T,X,Y,Z$ are $[X,Y]=Z$, $[T,X]=\frac{\alpha+2}{\alpha}X$, $[T,Y]=Y$ and $[T,Z]=2\frac{\alpha+1}{\alpha}Z$.
    
    Hence the Lie algebra generated by $X,Y,Z$ is isomorphic to $\heis$, and the Lie algebra $\lieg $ generated by $T,X,Y,Z$  is an extension of 
    $\heis$ by $\RR$.  Referring to the classification given in Section \ref{sec.heisenberg}, we see that:
    \begin{itemize}
     \item{If $\alpha=-1$, then $\lieg$ is isomorphic to $\RR \ltimes_s \heis$.}
     \item{If $\alpha \not = -1$, $\lieg$ is isomorphic to some algebra $\RR \ltimes_{\nu} \heis$.  Actually the parameter $\nu$  is given
      by the formula $\nu=\frac{det(A)}{Tr(A)^2}$, where $A$ is the matrix 
      $\left(  \begin{array}{cc} \frac{\alpha +2}{\alpha} & 0 \\ 0 & 1 \\  \end{array}\right)$. 
      Hence $\nu=\frac{1}{4}(1-\frac{1}{(\alpha+1)^2})$.  In particular, we see that the set of possible values for $\nu$ is exactly 
        $(- \infty, \frac{1}{4})$.}
    \end{itemize}

   The local isometric action of the Lie algebra $\lieg$ integrates into an action of a $4$-dimensional Lie group $G$ (we can directly describe $G$ as a subgroup
    of affine transformations of $\RR^3$) which is transitive on 
   the set $x_3>0$.  Hence  the metric $g_{\alpha}$ are always homogeneous.  Let us now 
   explain why the Lie algebra $\lieg$ coincides with the Lie algebra of local Killing fields of $g_{\alpha}$.  If this latter Lie algebra were bigger
    than $\lieg$, then it sould be $6$-dimensional (as already observed, there are no $3$-dimensional Lorentz metrics with a $5$-dimensional local Killing algebra),
     hence of constant curvature. It is however readily checked that $g_{\alpha}$ does not have constant curvature.
%
%

   Let now  $\nu_1$ and $\nu_2$ be as in the statement of Theorem \ref{theo.c1}.  We choose $\alpha_1$ and $\alpha_2$ in $(- \infty, -1)$
   such that  $\nu_1=\frac{1}{4}(1-\frac{1}{(\alpha_1 +1)^2}) $ and $\nu_2=\frac{1}{4}(1-\frac{1}{(\alpha_2 +1)^2}) $.
   There will be a unique $z \in (1,+ \infty)$ such that $\alpha_2=\alpha_1 \frac{z}{z-1}$.  Finally we introduce the metric
   $\tilde{g}_{\alpha_2}= z^{\alpha_1}(z-1)^{- \alpha_2} (x_3-1)^{\alpha_2}g_0$ (which is a smooth Lorentz metric for $x_3>1$). 
   
   On the open set $M$  given by $x_3>1$, we define a Lorentz metric $g$  as follows :
   \begin{itemize}
    \item{On the set $1<x_3\leq z$, $g=g_{\alpha_1}$.}
    \item{On the set $x_3>z$, $g=\tilde{g}_{\alpha_2}$.}
   \end{itemize}
We claim that $g$ has the properties of Theorem \ref{theo.c1}.   It is  smooth outside the set $x_3=z$.
The metrics $g_{\alpha_1}$ and $\tilde{g}_{\alpha_2}$ coincide along $x_3=z$.  Because of the condition $\alpha_2=\alpha_1 \frac{z}{z-1}$, their
 first derivatives also coincide at $x_3=z$, so that the metric $g$ is $C^1$. 
 On the open set $U_1$ defined by $1<x_3<z$, $g$ coincides with $g_{\alpha_1}$, hence its local Killing algebra is isomorphic to
  $\RR \ltimes_{\nu_1} \heis$.  On the open set $U_2$, the local Killing algebra of $g$ is that of $\tilde{g}_{\alpha_2}$.
   But if $\psi$  denotes the translation $x \mapsto x - e_3$, $\tilde{g}_{\alpha_2}$ is nothing but 
   $z^{\alpha_1}(z-1)^{-\alpha_2}\psi^*g_{\alpha_2}$, hence its local Killing algebra is isomorphic to that of $g_{\alpha_2}$, namely
    $\RR \ltimes_{\nu_2} \heis $.  This proves  Theorem \ref{theo.c1}.

    To get $C^0$ quasihomogeneous metrics on the $3$-torus, we do as follows.  We first choose $\alpha_1 \in (- \infty, -1) $ and $\alpha_2 \in 
    (1, + \infty)$ such that 
    $\nu_1=\frac{1}{4}(1-\frac{1}{(\alpha_1 +1)^2}) $ and $\nu_2=\frac{1}{4}(1-\frac{1}{(\alpha_2 +1)^2}) $ are not equal.  The
      two metrics $g_{\alpha_1}$ and $g_{\alpha_2}$ coincide on the hyperplane $x_3=1$.  Let us consider the strip $\overline{\Omega}$ given by
      $\frac{1}{2} \leq x_3 \leq 2^{\alpha_2 -\alpha_1}$, and define $g$ on $\overline{\Omega}$, by $g=g_{\alpha_1}$ if $\frac{1}{2} \leq x_3 \leq 1$
       and $g=g_{\alpha_2}$ on $1 \leq x_3 \leq 2^{\alpha_2-\alpha_1}$.  This yields a $C^0$ quasihomogeneous Lorentz metric on $\Omega$.  Now, observe
        that $g$ is invariant by the translations of vectors $e_1$  and $e_2$, and that the translation of vector 
        $(2^{\alpha_2-\alpha_1}-\frac{1}{2})e_3$ yields  an isometric identification of the restriction of $g$ to the hyperplanes
        $x_3=\frac{1}{2}$ and $x_3=2^{\alpha_2 - \alpha_1}$.  After those identifications, $g$ induces a $C^0$ quasihomogeneous metric on
        ${\mathbb T}^3$.  
        
        Observe that we could make more complicated examples, using  a countable family of $g_{\alpha_i}$, and producing 
         quasihomogenous metrics for which the locus of homogeneity has infinitely many connected components.  
  

\section{Annex A: Integrability locus}
\label{sec.annexeA}

Let $(M,\hm,\omega)$ be a Cartan geometry with  model space ${ X}=G/P$. We denote in the following by $m$ the dimension of $\lieg$.
If we fix $(e_1,\ldots,e_m)$  a basis of $\lieg$, we get a parallelism ${\mathcal P}$ on $\hm$   defined by the vector
fields $X_i=\omega^{-1}(e_i),\ i=1 \ldots,m $.  
As explained in \cite{vincent}[Sec. 4.3], it is enough to prove Theorem \ref{thm.integrabilite} for
the Cartan geometry  defined by $\mathcal P$ on $\hm$ (which is a Cartan geometry modelled on the abelian Lie group $\RR^m$).  Indeed, if we work on
 $\hm$, the Killing generators of order 
$r \geq 1$ are the same for both geometries, and the Killing fields are also the same. 

      \subsubsection{Curvature of the parallelism}
      Let $\mathcal P$ be the   parallelism on $\hm$ as defined above.  The curvature of ${\mathcal P}$
      is the $2$-form on $\hm$ defined by 
      $$K(X_i,X_j)=-\omega([X_i,X_j]), \text{for every pair of integers } 1 \leq i \leq j \leq n.$$  Let us write
      $[X_i,X_j]=\Sigma_{i,j} \gamma_{ij}^k X_k$, where the 
      $\gamma_{ij}^k$ are functions on  $\hm$.   At each $\hx \in \hm$, the curvature is seen as an element of $\Hom(\wedge^2(\lieg),\lieg)$
       given by the formula  $\kappa(\hx)(e_i,e_j)=\Sigma_{i,j}\gamma_{ij}^k(\hx)e_k$.  The curvature map 
        $\kappa : \hm \to \Hom(\wedge^2(\lieg),\lieg)$ of ${\mathcal P}$ differs from the curvature map of the initial Cartan geometry 
       $(M,\hm,\omega)$ by the constant bilinear form $[\, , \, ]_{\lieg}$, which explains that the derivatives
       (hence the Killing generators) of both curvature maps coincide.
       
       \subsubsection{Killing fields and  distribution on $\hm \times \RR^m$}
      
      The framing $(X_1,\ldots,X_m)$ identifies $T\hm$  with the product $\hm \times \RR^m$ in the following way: every vector 
      $\xi \in T_{\hx}\hm$ writes  $\xi_1 X_1(\hx)+\ldots+\xi_mX_m(\hx)$, allowing to identify  $\xi$  with  $(\xi_1,\ldots,\xi_m) \in \RR^m$.
      
      The local flows  $\varphi_{X_i}^t$ induce local flows on the tangent bundle $T\hm$, yielding vector fields $X_i^*$, $1 \leq i \leq m$, on 
   $\hm \times \RR^m$.  Those fields  have a simple expression involving $X_i$ and the curvature function 
   (computations can be found in \cite{vincent}[Lemma 4.9]) : 
   
    $$X_i^*(\hx,u)=(X_i(\hx),\kappa(\hx)(e_i,u)) \in T_{\hx}\hm \times \RR^m.$$
   
   Let us call $\Delta$ the distribution on $\hm \times \RR^m$, defined at each $(\hx,u)$ by
   $\Delta(\hx,u)=\operatorname{Span}(X_1^*(\hx,u), \ldots, X_m^*(\hx,u))$. The distribution $\Delta$ is related to Killing fields as follows.   The identification  $T\hm \simeq \hm \times \RR^m$ mentioned above 
     allows to see every local vector field on $U \subset \hm$ as a map $\varphi_X : U \to \RR^m$.  Assume
    that such a vector field satisfies $X(\hx)=u$. The condition that  $X$  commutes
      with the fields $X_i's$  up to order  $1$ at $\hx$, is equivalent to the  graph of the map $\varphi_X$ being tangent to $\Delta(\hx,u)$ at 
      $(\hx,u)$.   We thus see that if there exists an  integral manifold of $\Delta$ in the neighborhood of the point $(\hx,u)$, then
       this manifold is locally the graph of a map 
       $\varphi_X: U \subset \hm \to \RR^m$, associated to a Killing field $X$  on $U$, with $X(\hx)=u$  
     (see \cite{vincent}[Lemma 4.6]).   
       
       Observe  that except in the case of a {\it flat} parallelism (namely when all the $X_i$'s commute pairwise), the distribution 
       $\Delta$  never
         satisfies the Frobenius integrability condition  on $\hm \times \RR^m$.  The best we can hope is to find a few integral leaves for 
         $\Delta$.

%
%
      
      \subsubsection{Some useful submanifolds in $\hm \times \RR^m$}
     We recall from Sections \ref{sec.curvaturemap} and \ref{sec.killinggenerators} the definition of the derivatives $D^r \kappa$, the general
      curvature map $\dk$, and $\operatorname{Kill}^r(\hx)$, 
     the set of Killing generators of order $r$ at $\hx$.   Recall that $\hm^{\rm int}$ is the dense open subset of $\hm$ where the rank of
      $\dk$ is locally constant.  
   
   We define for every $r \in \NN^*$,
    $$ F_r:= \{ (\hx,u) \in \hm \times \RR^m \ | \  u \in \operatorname{Kill}^r(\hx) \}.$$
    Those are submanifolds of $\hm \times \RR^m$, as well as
    $$ \Sigma^{\rm int}:=( \hm^{\rm int} \times \RR^m) \cap F_{m+1}.$$
    
By what was said before,  proving Theorem \ref{thm.integrabilite} 
  is equivalent to finding an 
     integral leaf of $\Delta$ through each point of $\Sigma^{int}$. It is thus enough to show that $\Delta$ defines a distribution on
      $\Sigma^{\rm int}$ satisfying Frobenius integrability condition.  To do this, we first work on smaller manifolds, defining  
      $$\Omega:=\{ \hx \in \hm \ | \ \text{dim}(\operatorname{Kill}^j), \ 1\leq j\leq m+2, \  \text{are locally constant around $\hx$}  \},$$
   and for every $r \in \NN^*$,
   $$U_r:=\{ \hx \in \Omega \ | \ \dim(\operatorname{Kill}^r(\hx))=\dim(\operatorname{Kill}^{r+1}(\hx))  \},$$
  $$\Sigma_r=F_r \cap (U_r \times \RR^m).$$ 


%
 The first point is that $\Delta$ is tangent to the submanifolds $\Delta$, as shows the following lemma. 
 
 \begin{lemme}[\cite{nomizu}, Lemma 12; compare Lemma 4.7 of \cite{vincent}]
   For each $1 \leq r \leq m+1$, the submanifolds ${ \Sigma}_r$  are stable by the local flows $\phi_{X_i^*}^t$.  More precisely, 
   if $(\hx,u) \in { \Sigma}_r$, then for every $i=1,\ldots,n$,
   $X_i^*(x,u) \in T {\Sigma}_r$.
   \end{lemme}
   
 
 Next, because $F_r \subset F_1$, the distribution $\Delta$, seen as a distribution of $\Sigma_r$ satisfies the Frobenius condition,
    hence is integrable. This follows from
   \begin{lemme}[\cite{vincent}, Lemma 4.10] 
   \label{lem.involutif}
   At each $(\hx,u) \in F_1$ and for every $1 \leq i \leq j \leq m$, the brackets $[X_i^*(x,u),X_j^*(x,u)]$   belong to $\Delta(x,u)$.
   \end{lemme}
   
    From this integrability property, and the previous discussion,  we infer that  $\dim \operatorname{Kill}^r(\hx) \geq \dim \kiloc(\hx)$ 
    for all $\hx \in U_r$.  Since the reverse inequality always holds, we get for every $1 \leq r \leq m+1$ :
    $$\operatorname{Kill}^{r}(\hx)= \operatorname{Kill}^{m+1}(\hx)=\kiloc(\hx).$$ 
     We conclude that for every  $1 \leq r \leq m+1$ ,  $ \Sigma^{int} \cap (U_r \times \RR^m)$ is an open subset of $ \Sigma_r $. 
       
       We now conclude in the following way. The inequalities $\dim (\operatorname{Kill}^j)\geq \dim (\operatorname{Kill}^{j+1})$ show
        that at each point $\hx \in \Omega$, there must exist
        $1 \leq r \leq m+1$ such that $\dim (\operatorname{Kill}^r)=\dim (\operatorname{Kill}^{r+1})$. This   implies that
        $ \Sigma^{int} \cap (\Omega \times \RR^m)$ is an open subset of $\bigcup_{r=1}^{m+1}\Sigma_r $.  In particular, $\Delta$ 
        is tangent to $\Sigma^{int} \cap (\Omega \times \RR^m)$ and satisfies the Frobenius integrability condition there.
         Finally, because $\Sigma^{int} \cap (\Omega \times \RR^m)$ is a dense open set of the manifold $\Sigma^{int}$, this property holds true
          on $\Sigma^{int}$.    This yields Theorem \ref{thm.integrabilite}.
         


    \section{Annex B: Closed orbits in representations of $\PSL(2,\RR)$}
    \label{annexeB}
 
We provide below a proof of Lemma \ref{lem.orbitehyper}.

Let us consider the  basis of $\sld$ given by  $$H=\left( \begin{array}{cc} 1&0 \\ 0& -1  \end{array} \right), 
E=\left( \begin{array}{cc}0 &1 \\ 0& 0  \end{array} \right) \text{ and } F=\left( \begin{array}{cc} 0&0 \\ 1& 0  \end{array} \right).$$  
All hyperbolic flows in 
 $\PSL(2,\RR)$ are conjugated inside $\PSL(2,\RR)$, and the same is true for elliptic flows.  Hence, we may assume that $\{ \varphi^t \}$ is generated either by $H$, or by 
  $E-F$.

  We first do the proof when the representation is irreducible.  We then know that up to isomorphism, $\rho$ is induced by the action of $\SL(2,\RR)$ on homogeneous
 polynomials in $2$ variables, the action being by linear substitution. Observe that the action of $-Id \in \PSL(2,\RR)$ is trivial only when the polynomials
 have even degree, so that $V$ is odd dimensional.  Let us put $2k+1=\rm{dim } V$, and denote by 
 $\overline{\rho} : \sld \to \mathfrak{gl}(2k+1,\RR)$  the representation induced by $\rho$ at the Lie algebra level. There
 is a suitable basis $e_1,\ldots,e_{2k+1}$ of $V$ where 
  
$$\overline{H}=\overline{\rho}(H)=\left( \begin{array}{ccccc} 2k & & & & \\
 & 2k-2 & & & \\ & & \ddots & & \\
  & & & -2k+2 & \\
   & & & &-2k \end{array} \right) $$
   
   $$\overline{E}-\overline{F}=\overline{\rho}(E-F)=\left( \begin{array}{ccccc} 0 & 1 &0 & & \\
 -2k & \ddots & 2&\ddots & \\ &-2k+1 & \ddots & \ddots&0 \\
  & & \ddots& \ddots & 2k\\
   & & &-1 &0 \end{array} \right) $$
   

Because $2k+1$ is odd, it is easy to check that $\rm{Fix}(\varphi^t)$ consists of  the line generated by $v$.  Moreover, 
$\rho(\SL(2,\RR))$ preserves a pseudo-Riemannian scalar product $g$ on $V$, having type $(k,k+1)$.  Precisely, 
$ g(x,x)=(\Sigma_{m=1}^k 2a_mx_mx_{2k+2-m})+a_{k+1}x_{k+1}^2$, where $a_1$ is any element of $\RR^*$ and the $a_i$'s satisfies the relations
 $a_{i+1}=\frac{2k+1-i}{i}a_i$  for $i=1,\ldots,k$.  
 
 Let us check that $g(v,v) \not =0$.  When $\{\varphi^t\}$ is generated by $H$, it is obvious since   $v$ is then on the line $\RR.e_{k+1}$.  
  Assume now that $\varphi^t$ is elliptic.  If $g(v,v)=0$, then $v$ is included in the $g$-orthogonal  $ {\{v \}}^{\perp}$.
  By compacity of the $1$-parameter group $\{ \varphi^t\}$, we
    get a $\varphi^t$-invariant 
  decomposition   $v^{\perp}=L\oplus \RR v$.  Now, the dimension of $L$ is odd (namely $2k-1$),   
  hence $\varphi^t$  admits a line of fixed points in $L$. This gives a $2$-plane of fixed points for $\varphi^t$ in $V$: contradiction.
  
  We can conclude that the orbit $\O$ of $v$ under $\PSL(2,\RR)$ is closed in the following way.  Because the action of $\SL(2,\RR)$ is algebraic, 
   a non closed orbit must accumulate on some orbit of (strictly) smaller dimension.  In our case, if not closed, the orbit $\O$ should have the origin in 
   its closure,
     because there are no orbits of dimension $1$ for finite dimensional representations of $\PSL(2,\RR)$. This is not possible, 
      since  $\mathcal{O}$ is included in the set $\{ x \in V \ | \ g(x,x)=g(v,v)  \}$ and $g(v,v) \not =0$. 

  In the  case $V$ is not irreducible,it splits as a direct sum $V=V_1 \oplus \ldots \oplus V_s$ of irreducible representations.
 Let $u \in V$ having stabizer $\{ \varphi^t \}$, a hyperbolic or an elliptic flow.  We write $u=u_1+\ldots+u_s$.   
 In any finite dimensional irreducible representation of $\PSL(2,\RR)$, the stabilizer of a point is either
   $\PSL(2,\RR)$, a $1$-parameter group, or $\{ id \}$.  Hence, the stabilizer of each $u_i$ is either
  $\PSL(2,\RR)$ or $\{ \varphi^t \}$, because if other possibilities occured, the stabilizer of $u$ would be trivial.  We assume this stabilizer
  is 
  $\PSL(2,\RR)$ for $u_1,\ldots,u_l$ and $\{ \varphi^t \}$ for $u_{l+1}, \ldots, u_s$.  Let $(g_k)$ be a sequence in $\PSL(2,\RR)$  such that
   $g_k.u$ converges to $u_{\infty} \in V$.  For each $i=l+1, \ldots ,s$, we have $g_k.u_i \to u_i^{\infty}$. Because $V_i$ is irreducible, we know by
    the first part of the proof that the orbit of $u_i$ under 
   $\PSL(2,\RR)$ is closed.  Hence there exists $\tilde{g}_k$ converging to $g_{\infty} \in \PSL(2,\RR)$ such that 
   $\tilde{g}_k.u_{l+1}=g_k.u_{l+1}$.  In particular $g_k=\tilde{g}_k.\varphi^{t_k}$ for all $k$.  But now, this implies
   $g_k.u_i=\tilde{g}_k.u_i$
    for all $i=l+1, \ldots, s$, hence for all $i=1,\ldots,s$.  We end up with $g_k.u=\tilde{g}_k.u$, so that $u_{\infty}=g_{\infty}.u$.  
    This concludes the proof.

\end{document}